\documentclass[11pt]{amsart}

\usepackage[all]{xypic}

\usepackage{latexsym}

\usepackage{amssymb}
\usepackage{amsfonts}
\usepackage{mathrsfs}
\usepackage{eucal}
\usepackage{amscd}
\usepackage{amsmath,amsthm}

\newtheorem{lemma}{Lemma}[section]

\newtheorem{lem}[lemma]{Lemma}
\newtheorem{prop}[lemma]{Proposition}
\newtheorem{thm}[lemma]{Theorem}
\newtheorem{cor}[lemma]{Corollary}

{

}

\theoremstyle{definition}

\theoremstyle{remark}

\numberwithin{equation}{section}

\newenvironment{pf}{\noindent{\bf Proof.}}{\hfill $\square$\medskip}

\def\CC{{\mathbb C}}

\def\NN{{\mathbb N}}

\def\PP{{\mathbb P}}

\def\ZZ{{\mathbb Z}}

\def\Mol{{\bar M}}

\def\0ol{{\bar 0}}
\def\1ol{{\bar 1}}
\def\2ol{{\bar 2}}
\def\ol2{{\bar 2}}
\def\3ol{{\bar 3}}
\def\4ol{{\bar 4}}
\def\5ol{{\bar 5}}
\def\6ol{{\bar 6}}
\def\7ol{{\bar 7}}
\def\8ol{{\bar 8}}
\def\9ol{{\bar 9}}

\def\bold0{{\bf 0}}
\def\bold1{{\bf 1}}
\def\bold2{{\bf 2}} 
\def\bold3{{\bf  3}}
\def\bold4{{\bf 4}}
\def\bold5{{\bf 5}}
\def\bold6{{\bf 6}}
\def\bold7{{\bf 7}}
\def\bold8{{\bf 8}}
\def\bold9{{\bf 9}}

\def\P2Skly{\PP^2_{Skly}}

\def\add{  {\sf add}}

\def\End{\operatorname {End}}

\def\gr{\operatorname {gr}}

\def\Hom{\operatorname {Hom}}

\def\ker{\operatorname {ker}}

\def\op{{\operatorname {op}}}

\def\th{\operatorname {th}}    

\def\dim{\operatorname{dim}}

\def\End{\operatorname{End}}

\def\Fdim{{\sf Fdim}}
\def\fdim{{\sf fdim}}

\def\Gr{{\sf Gr}}

\def\gr{{\sf gr}}

\def\Hom{\operatorname{Hom}}

\def\id{\operatorname{id}}

\def\liminj{\varinjlim}

\def\mod{{\sf mod}}
\def\Mod{{\sf Mod}}

\def\op{\operatorname{op}}

\def\Proj{\operatorname{Proj}}

\def\proj{\operatorname{\sf proj}}

\def\QGr{\operatorname{\sf QGr}}
\def\qgr{\operatorname{\sf qgr}}

\def\Span{\operatorname{span}}

\def\ul1{\operatorname{\underline{1}}}

\def\l{\leftarrow}

\def\d{\downarrow}

\def\a{\alpha}

\def\d{\delta}

\def\l{\lambda}

\def\s{\sigma}

\def\L{\Lambda}

\def\sA{{\sf A}}

\def\sD{{\sf D}}

\def\sK{{\sf K}}

\def\sT{{\sf T}}

\def\cal{\mathcal}

\def\cF{{\cal F}}

\def\cH{{\cal H}}

\def\cM{{\cal M}}

\def\cO{{\cal O}}
\def\cP{{\cal P}}

\def\cS{{\cal S}}

\def\Qcoh{{\sf Qcoh}}

\def\sRings{{\sf Ring}}

\def\dirlim{\mathop{\vtop{\baselineskip -100pt\lineskip -1pt\lineskiplimit 0pt
\setbox0\hbox{lim}\copy0\hbox to \wd0{\rightarrowfill}}}\limits}
\def\invlim{\mathop{\vtop{\baselineskip -100pt\lineskip -1pt\lineskiplimit 0pt
\setbox0\hbox{lim}\copy0\hbox to \wd0{\leftarrowfill}}}\limits}

\def\I11{{1 \kern -0.8pt \! \mbox{l}}}
\def\mumu{{\mu\kern-4.2pt\mu}}
\def\bfmu{{\mu\kern-4.2pt\mu}}
\def\2slash{\backslash \! \backslash}

\def\boxtimes{\setbox0\hbox{$\Box$}\copy0\kern-\wd0\hbox{$\times$}}

\def\hdot{{\:\raisebox{2pt}{\text{\circle*{1.5}}}}}

\date{}                                           

\begin{document}

\title[Category equivalences involving path algebras of quivers]{Category equivalences involving  graded modules over path algebras of quivers}

\author{S. Paul Smith}

\address{ Department of Mathematics, Box 354350, Univ.
Washington, Seattle, WA 98195}

\email{smith@math.washington.edu}

\thanks{This work was partially supported by the U.S. National Science Foundation grant 0602347. }


\begin{abstract}
Let $Q$ be a finite quiver with vertex set $I$ and arrow set $Q_1$, $k$ a field, and $kQ$ its path algebra with its standard grading. 
This paper proves some category equivalences involving the quotient category $\QGr(kQ) := \Gr(kQ)/\Fdim(kQ)$ of graded $kQ$-modules modulo those that are the sum of their finite dimensional submodules, namely
$$
 \QGr(kQ) \equiv \Mod S(Q) \equiv \Gr L(Q^\circ) \equiv \Mod L(Q^\circ)_0 \equiv \QGr(kQ^{(n)}).
 $$
 Here $S(Q)=\liminj \End_{kI}(kQ_1^{\otimes n})$ is a direct limit of finite dimensional semisimple algebras; $Q^\circ$ is the quiver without sources or sinks that is obtained by repeatedly removing all sinks and sources from $Q$; $L(Q^{\circ})$ is the Leavitt path algebra of $Q^{\circ}$; $L(Q^\circ)_0$ is its degree zero component; 
 and $Q^{(n)}$ is the quiver whose incidence matrix is the $n^{\th}$ power of that for $Q$. 
 It is also shown that all short exact sequences in $\qgr(kQ)$, the full subcategory of finitely presented objects in $\QGr(kQ)$, split. Consequently $\qgr(kQ)$ can be given the structure of a triangulated category with suspension
 functor the Serre degree twist $(-1)$; this triangulated category is equivalent to the ``singularity category''
 $\sD^b(\L)/\sD^{\sf perf}(\L)$ where $\L$ is the radical square zero algebra $kQ/kQ_{\ge 2}$, and $\sD^b(\L)$ is the bounded derived category of finite dimensional left $\L$-modules.
\end{abstract}

\maketitle

\pagenumbering{arabic}

\setcounter{section}{0}

\section{Introduction}

\subsection{}
\label{sect.1.1}
Throughout $k$ is a field and
$Q$ a finite quiver (directed graph) with vertex set $I$. Loops and multiple arrows between vertices are allowed. 

We write $kQ$ for the path algebra of $Q$. 

We make $kQ$ an $\NN$-graded algebra by declaring that a path is homogeneous of degree equal to its length.  The category of $\ZZ$-graded left $kQ$-modules with degree-preserving homomorphisms is denoted by $\Gr(kQ)$ and we write $\Fdim( kQ)$ for its full subcategory of modules that are the sum of their
finite-dimensional submodules. Since $\Fdim( kQ)$ is a localizing subcategory of $\Gr(kQ)$ we may form the 
quotient category 
$$
\QGr (kQ) := \frac{\Gr (kQ) }{\Fdim  (kQ) }.
$$
By \cite[Prop. 4, p. 372]{Gab}, the quotient functor $\pi^*:\Gr (kQ)  \to \QGr (kQ)$ has a right adjoint  
that we will denote by $\pi_*$. We define 
$$
\cO:=\pi^*(kQ).
$$

The main result in this paper is the following theorem combined with an explicit description of the
algebra $S(Q)$ that appears in its statement. 

\begin{thm}
\label{thm1}
The endomorphism ring of $\cO$ in $\QGr(kQ)$ is an ultramatricial $k$-algebra, 
$S(Q)$, and $\Hom_{\QGr(kQ)}(\cO,-)$ is an equivalence  
$$
 \QGr (kQ)\equiv \Mod S(Q)
 $$ 
 with the category of right $S(Q)$-modules. 
\end{thm}

\subsection{Definition and description of $S(Q)$}
\label{sect.defn.SQ}

We write $Q_n$ for the set of paths of length $n$ and 
$kQ_n$ for the linear span of $Q_n$. With this notation
\begin{align*}
kQ = & kI \oplus kQ_1 \oplus kQ_2 \oplus \cdots
\\
= &  T_{kI}(kQ_1)
\end{align*}
where $ T_{kI}(kQ_1)$ is the tensor algebra of the $kI$-bimodule $kQ_1$. 

The ring of left $kI$-module endomorphisms of $kQ_n$ is denoted
$$
S_n := \End_{kI} (kQ_n).
$$  
Since $kQ_{n+1}=kQ_1 \otimes_{kI} kQ_n$
the functor $kQ_1 \otimes_{kI}-$ gives $k$-algebra homomorphisms $$\theta_n:S_n \to S_{n+1}.$$ 
Explicitly, if $x_1,\ldots,x_{n+1} \in kQ_1$, $f \in S_n$, and $\otimes=\otimes_{kI}$, then
$$
\theta_n(f)\big(x_1 \otimes \cdots \otimes x_{n+1}\big) := 
x_1 \otimes f(x_2  \otimes \cdots \otimes x_{n+1}).
$$
The $\theta_n$s give rise to a directed system $kI = S_0 \to S_1 \to \cdots$, and we define
$$
S(Q): = \liminj S_n.
$$

As a $k$-algebra, $kI$ is isomorphic to a product of $|I|$ copies of $k$. Every left $kI$-module is therefore a direct sum of 1-dimensional $kI$-submodules and the endomorphism ring of a finite dimensional $kI$-module
is therefore a direct sum of $\le |I|$ matrix algebras $M_r(k)$ where the $r$s that appear are determined by
the multiplicities of the simple $kI$-modules.

Hence $S(Q)$ is a direct limit of products of matrix algebras. Such   algebras are called {\sf ultramatricial}. 

Theorem \ref{thm1} will follow from the following result.

\begin{thm}
\label{thm2}
The object $\cO$ is a finitely generated projective generator in $\QGr (kQ) $ and 
$$
\End_{\QGr(kQ) } \cO \cong S(Q).
$$
The functor implementing the equivalence in Theorem \ref{thm1} is $\Hom_{\QGr(kQ) }(\cO,-)$. 
\end{thm}
 
Ultramatricial algebras are described by Bratteli diagrams \cite{B} (see also \cite{KG2}). 
The Bratteli diagram for $S(Q)$, and hence $S(Q)$, is described explicitly in Proposition  \ref{prop.brat} 
 in terms of the incidence matrix for $Q$.

\subsection{Relation to Leavitt path algebras and Cuntz-Krieger algebras}

Apart from the path algebra $kQ$ two other algebras are commonly associated to $Q$, the Leavitt
path algebra $L(Q)$ and the Cuntz-Krieger algebra ${\mathscr O}_Q$.  The algebra $L(Q)$, which can be defined over any commutative ring, is an algebraic analogue of  the C$^*$-algebra ${\mathscr O}_Q$ 
because (when the base field is $\CC$)  ${\mathscr O}_Q$ contains $L(Q)$ as a dense subalgebra.

\begin{thm}
\label{thm.LPAs}
Let $Q^\circ$ be the quiver without sources or sinks that is obtained by repeatedly removing all sinks and sources from $Q$. Then
\begin{enumerate}
  \item 
  $\QGr(kQ) \equiv \QGr(kQ^\circ)$;
  \item 
 $S(Q^\circ) \cong L(Q^\circ)_0$; 
  \item 
 $L(Q^\circ)$ is a strongly graded ring;
 \item{}
   $\QGr(kQ) \equiv \Mod S(Q) \equiv \Gr L(Q^\circ) \equiv \Mod L(Q^\circ)_0$.
\end{enumerate}
\end{thm}

 After proving this theorem the author learned that Roozbeh Hazrat had previously given 
 necessary and sufficient conditions for  $L(Q^\circ)$ to be a strongly graded ring  \cite[Thm. 3.15]{H}. 
 The idea in our proof of (3) differs from that in Hazrat's paper.

\subsection{Coherence}
A ring $R$ is {\sf left coherent} if the kernel of every homomorphism $f:R^m \to R^n$ between finitely generated 
free left $R$-modules is finitely generated. 
If $R$ is left coherent we write $\mod R$ for the full subcategory of $\Mod R$ consisting of 
finitely presented modules; $\mod R$ is then an abelian category.

To prove $R$ is left coherent it suffices to show that every finitely generated left ideal is finitely presented. 

A ring in which every left ideal is projective is left coherent so $kQ$ is left coherent. 
A direct limit of left coherent rings is left coherent so $S(Q)$ is left coherent.

Because $kQ$ is left coherent the full subcategory
$$
\gr(kQ) \; \subset \; \Gr(kQ)
$$
consisting of finitely presented graded left $kQ$-modules is abelian.   
The category $\fdim \big( kQ\big):=\big(\gr(kQ) \big) \cap \big(\Fdim (kQ) \big)$ is the full subcategory of $\gr(kQ)$ consisting 
of finite dimensional modules. We now define 
define 
$$
\qgr (kQ ) := \frac{\gr(kQ)}{\fdim (kQ)} \; \; \subset \QGr (kQ).
$$
 
\begin{prop}
The equivalence in Theorem \ref{thm1} restricts to an equivalence
$$
\qgr (kQ) \equiv \mod S(Q).
$$
\end{prop} 
 
By \cite[Prop. A.5, p. 113]{HK}, $\qgr R$ consists of the finitely presented objects in $\QGr R$
and every object in  $\QGr R$ is a direct limit of objects in $\qgr R$.\footnote{An object $\cM$ in an additive category $\sA$ is finitely presented if $\Hom_\sA(\cM,-)$ commutes with direct limits; is finitely generated if whenever $\cM = \sum \cM_i$ for some directed family of subobjects $\cM_i$ there is an index $j$ such that $\cM=\cM_j$; is coherent if it is finitely presented and all its finitely generated subobjects are finitely presented.}

\subsection{$\qgr (kQ)$ as a triangulated category}

One of the main steps in proving Theorems \ref{thm1} and \ref{thm2} is to prove the following.

\begin{prop}
Every short exact sequence in $\qgr(kQ)$ splits.
\end{prop}

If $\Sigma$ is an auto-equivalence of an abelian category $\sA$  in which every short exact sequence splits,
then $\sA$ can be given the structure of a triangulated category with $\Sigma$ being the translation: 
one declares that the distinguished triangles are all direct sums of the following triangles:
\begin{align*}
& \cM \to 0 \to \Sigma \cM \stackrel{\id}{\longrightarrow} \Sigma \cM, 
\\
& \cM \stackrel{\id}{\longrightarrow} \cM  \to 0 \to \Sigma \cM,
\\
& 0 \to \cM \stackrel{\id}{\longrightarrow} \cM  \to 0,
\end{align*}
as $\cM$ ranges over the objects of $\sA$.
Hence $\qgr(kQ)$ endowed with the Serre twist $(-1)$ is a triangulated category.  In section \ref{sect.chen} 
we use a result of Xiao-Wu Chen to prove the following.

\begin{thm}
Let $Q$ be a quiver and $\L$ the finite dimensional algebra $kQ/kQ_{\ge 2}$. 
There is an equivalence of triangulated categories
$$
\big(\qgr(kQ), (-1)\big) \equiv \sD^b(\mod \L)/\sD^b_{\rm perf}(\mod \L).
$$
\end{thm}

\subsection{Equivalences of categories}

It can happen that $\QGr(kQ)$ is equivalent to $\QGr (kQ')$ with $Q$ and $Q'$ being non-isomorphic quivers. 

\begin{thm}
\label{thm3}
(See section \ref{sect.sink.source}.)
If $Q$ and $Q'$ become the same after repeatedly removing vertices that are sources or sinks, then $\QGr(kQ)  \equiv \QGr(kQ')$.
\end{thm}

Let $A$ be a $\ZZ$-graded algebra. If $m$ is a positive integer the algebra $A^{(m)}=\oplus_{i \in \ZZ} A_{im}$
is called the $m^{\th}$ Veronese subalgebra of $A$. When $A$ is a commutative $\NN$-graded algebra the schemes $\Proj A$ and $\Proj A^{(m)}$ are isomorphic.
Verevkin  proved a non-commutative version of that result: $\QGr A \equiv \QGr A^{(m)}$ if $A$ is an $\NN$-graded ring  generated by $A_0$ and $A_1$  \cite[Thm. 4.4]{V}.

\begin{thm}
\label{thm4}
Let $Q$ be a quiver with incidence matrix $C$. Let $Q^{(m)}$ be the quiver with incidence matrix $C^m$, $m \ge 1$; i.e., $Q^{(m)}$ has the same vertices as $Q$ but the arrows in $Q^{(m)}$ are the paths of length $m$ in $Q$.
Then $$\QGr(kQ)  \equiv \QGr(kQ^{(m)}).$$
\end{thm}
\begin{pf}
This follows from Verevkin's result because $kQ^{(m)} = (kQ)^{(m)}$. 
\end{pf}

The referee pointed out the following alternative proof of Theorem \ref{thm4}. 
First, $S_n(Q^{(m)})=S_{nm}(Q)$ so the directed system used to define $S(Q^{(m)})$ is equal to the
directed system obtained by taking every $m^{\th}$ term of the directed system for $S(Q)$. Hence 
$S(Q) = S(Q^{(m)})$. Therefore
$$
\QGr(kQ) \equiv \Mod S(Q) = \Mod  S(Q^{(m)}) \equiv \QGr(kQ^{(m)}).
$$

We call $Q^{(n)}$ the {\sf $n^{\th}$ Veronese} of $Q$. In symbolic dynamics $Q^{(n)}$ is called the {\sf $n^{\th}$ higher power graph} of $Q$ \cite[Defn. 2.3.10]{LM}.

Other equivalences involving strong shift equivalence of incidence 
matrices, a notion from symbolic dynamics, appear in \cite{Sm3}.

\medskip
{\bf Acknowledgements.} 
Part of this paper was written during a visit to Bielefeld University in April 2011. I thank Henning Krause for the invitation and the university for providing  excellent working conditions.  

Conversations with Gene Abrams, Pere Ara, Ken Goodearl, Henning Krause, Mark Tomforde, and Michel Van den Bergh, shed light on the material in this paper and I thank them all for sharing their ideas and knowledge. 

I am especially grateful to Xiao-Wu Chen, Roozbeh Hazrat, and the referee, for reading earlier versions of this paper and pointing out typos and obscurities. Their comments and suggestions have 
improved the final version of this paper. I thank Chen for sending me an early version of his paper \cite{XWChen} and Hazrat for bringing his 
paper \cite{H} to my attention (see the remark after Theorem \ref{thm.LPAs}).


\section{The endomorphism ring of $\cO$}

Recall that $\cO$ denotes $\pi^*( kQ)$, the image of the graded left module $kQ$ in the quotient category
$\QGr(kQ)$. 

\smallskip

 {\bf Notation.}
 In addition to the notation set out at the beginning of section \ref{sect.defn.SQ} we write $Q_{\ge n}$ for the set 
 of paths of length $\ge n$ and  $kQ_{\ge n}$  for its linear span. We note that $kQ_{\ge n}$ is a
 graded two-sided ideal in $kQ$.  
 
 We write $e_i$ for the trivial path at vertex $i$, $E_i$ for the simple module at vertex $i$, and 
 $P_i=(kQ)e_i$. 
 
 If $p$ is a path in $Q$ we write $s(p)$ for its starting point and $t(p)$ for the vertex at which it terminates. 

We write $pq$ to denote the path $q$ followed by the path $p$.

\begin{lem}
\label{lem.tors.kQ} 
Let $I^0=\{i \in I \; | \; \hbox{the number of paths starting at $i$ is finite}\}$. Let $I^\infty:=I-I^0$ and let $Q^\infty$ be the subquiver of $Q$ consisting of the vertices in $I^\infty$ and all arrows that begin and end at points in
 $I^\infty$.
Let $T$ be the sum of all finite-dimensional left ideals in $kQ$. Then 
\begin{enumerate}
  \item 
  $T$ is a two sided ideal; 
  \item 
  $T=(e_i \; | \; i \in I^0)$;
  \item 
  $kQ/T \cong kQ^\infty$;
  \item{}
  the only finite-dimensional left ideal in $kQ^\infty$ is $\{0\}$;
  \item{}
  if $f:kQ_{\ge n} \to T$ is a homomorphism of graded left $kQ$-modules, then $f(kQ_{\ge n+r})=0$
  for $r \gg 0$. 
\end{enumerate}
\end{lem}
\begin{pf}
The result is obviously true if $\dim_k kQ <\infty$ so we assume this is not the case, i.e., $Q$ has arbitrarily
long paths; equivalently,   $Q^\infty \ne \varnothing$. 

(1)
If $L$ is a finite-dimensional left ideal in $kQ$ so is $Lx$ for all $x \in kQ$, whence $T$ is a two-sided ideal.

(2), (3), (4). 
Since the paths beginning at  a vertex $i$ are a basis for $(kQ)e_i$, $\dim_k (kQ)e_i <\infty$ if and only if $i \in I^0$. Hence $T$ contains $\{ e_i \; | \; i \in I^0\}$. 
It is clear that 
$$
\frac{kQ}{( e_i \; | \; i \in I^0)} \cong kQ^\infty.
$$

Let $p$ be a path in $Q^\infty$. Then there is an arrow $a \in Q^\infty$ such that $ap \ne 0$. It follows that 
$\dim_k(kQ^\infty)p =\infty$. It follows that the  only finite-dimensional left ideal in $kQ^\infty$ is $\{0\}$. Therefore $T/(e_i \; | \; i \in I^0)=0$. 

(5)
Let $f:kQ_{\ge n} \to T$ be a homomorphism of graded left $kQ$-modules.
Every finitely generated left ideal contained in $T$ has finite dimension so $f(kQ_{\ge n})$ has finite dimension. Hence $kQ_{\ge n}/\ker f$ is annihilated by $kQ_{\ge r}$ for $r \gg 0$. In other words, $\ker f \supset (kQ_{\ge r}).(kQ_{\ge n})=kQ_{\ge n+r}$. 
\end{pf}

The ideal $T$ in Lemma \ref{lem.tors.kQ} need not have finite dimension; for example, if $Q$
is the quiver in Proposition \ref{eg.sink.source}, $\dim_kT=\infty$.  

\begin{lemma}
\label{lem.cofinal}
Let $I$ be a graded left ideal of $kQ$. 
If $kQ/I$ is the sum of its finite dimensional submodules, then $I\supset kQ_{\ge n}$ for $n \gg 0$. 
\end{lemma}
\begin{pf}
The image of 1 in $kQ/I$ belongs to a {\it finite} sum of finite dimensional submodules of $kQ/I$ so
the submodule it generates is finite dimensional. Hence $\dim_k (kQ/I)<\infty$. Therefore $kQ/I$ is non-zero in only finitely many degrees; thus $I$ contains $kQ_{\ge n}$ for $n \gg 0$.  
\end{pf}

By definition, the objects in $\QGr(kQ)$ are the same as those in the $\Gr(kQ)$ and the morphisms are 
$$
\Hom_{\QGr(kQ)} (\pi^* M,\pi^* N) = \liminj \Hom_{\Gr(kQ)} (M',N/N')
$$
where the direct limit is taken as $M'$ and $N'$ range over all graded 
submodules of $M$ and $N$  
such that $M/M'$ and $N'$ belong to $\Fdim(kQ)$.

\begin{prop}
If $N \in \Gr(kQ)$, then
$$
\Hom_{\QGr(kQ)}(\cO,\pi^* N) = \liminj \Hom_{\Gr(kQ)} (kQ_{\ge n},N/N')
$$
where the direct limit is taken over all integers $n \ge 0$ and all submodules $N'$ of $N$  
such that $N'$ is the sum of its finite dimensional submodules.
\end{prop}
\begin{pf}
This follows from Lemma \ref{lem.cofinal}.
\end{pf}

   \begin{lem}
   \label{lem.2.4}
   Consider $kQ_n$ as a left $kI$-module.
The restriction map
$$
\Phi:\End_{\Gr(kQ) }(kQ_{\ge n}) \longrightarrow \End_{kI} (kQ_n), \qquad \Phi(f)=f|_{kQ_n},
$$
is a $k$-algebra isomorphism with inverse given by 
applying the functor $kQ \otimes_{kI} -$ to each $kI$-module endomorphism of $kQ_n$. 
\end{lem}
\begin{pf}
Each $f \in \End_{\Gr(kQ)}(kQ_{\ge n})$  sends $kQ_n$ to itself. Since $f$ is a left $kQ$-module homomorphism it is a left $kI$-module homomorphism. Hence $\Phi$ is a well-defined algebra homomorphism. Since $kQ_{\ge n}$ is generated by $kQ_n$ as a left $kQ$-module, $\Phi$ is
injective. Since $kQ_{\ge n} \cong kQ \otimes_{kI} kQ_n$ every left $kI$-module homomorphism 
$kQ_n \to kQ_n$ extends in a unique way to a $kQ$-module homomorphism 
$kQ_{\ge n} \to kQ_{\ge n}$ (by applying the functor $kQ \otimes_{kI} -$). Hence $\Phi$ is surjective.
\end{pf}
 
 \begin{thm} 
 \label{thm.End.O}
 There is a $k$-algebra isomorphism
  $$
 \End_{\QGr(kQ)} \cO \cong \liminj \Hom_{\Gr(kQ)}(kQ_{\ge n},kQ_{\ge n}) = S(Q).
 $$
 \end{thm}
 \begin{pf}
 By the definition of  morphisms in a quotient category,
\begin{equation}
\label{.d.sys1}
\End_{\QGr(kQ) }\cO  = \liminj \Hom_{\Gr(kQ)}(I,kQ/T')
\end{equation}
where $I$ runs over all graded left ideals such that $\dim_k(kQ/I) < \infty$ and 
$T'$ runs over all graded left ideals such that $\dim_k T' < \infty$.

If $T'$ is a graded left ideal of finite dimension it is contained in the ideal $T$ 
that appears in Lemma \ref{lem.tors.kQ}.  
The system of graded left ideals of finite codimension in $kQ$
 is cofinal with the system of left ideals
$kQ_{\ge n}$.  These two facts imply that  
$$
\End_{\QGr(kQ)} \cO = \liminj_{n} \Hom_{\Gr(kQ)}(kQ_{\ge n},kQ/T).
$$

 Since $kQ_{\ge n}$ is a projective left $kQ$-module the map
 $$
 \Hom_{\Gr(kQ)}(kQ_{\ge n},kQ) \to \Hom_{\Gr(kQ)}(kQ_{\ge n},kQ/T)
 $$
 is surjective. This leads to a surjective map
 \begin{equation}
 \label{eq.2.dlims}
 \liminj_{n} \Hom_{\Gr(kQ)}(kQ_{\ge n},kQ)  \to \liminj_{n} \Hom_{\Gr(kQ)}(kQ_{\ge n},kQ/T).
 \end{equation}
 Suppose the image of a map $f \in  \Hom_{\Gr(kQ)}(kQ_{\ge n},kQ)$ is contained in $T$.
 By Lemma \ref{lem.tors.kQ}(5), the restriction of $f$ to $kQ_{\ge n+r}$ is zero for $r \gg 0$. The map in 
 (\ref{eq.2.dlims}) is therefore injective and hence an isomorphism. 
 
 Since morphisms in $\Hom_{\Gr(kQ)}(kQ_{\ge n},kQ)$ preserve degree the natural map 
 $\Hom_{\Gr(kQ)}(kQ_{\ge n},kQ_{\ge n}) \to \Hom_{\Gr(kQ)}(kQ_{\ge n},kQ)$ is an isomorphism. It follows that 
$$
\End_{\QGr(kQ)} \cO = \liminj_{n} \Hom_{\Gr(kQ)}(kQ_{\ge n},kQ_{\ge n}).
$$
However, $\Hom_{\Gr(kQ)}(kQ_{\ge n},kQ_{\ge n}) \cong \End_{kI}(kQ_n)$ by Lemma \ref{lem.2.4} so the result follows from the definition of $S(Q)$. 
 \end{pf}
 
 
\section{Proof that $\cO$ is a progenerator in $\QGr(kQ)$}
\label{sect.O}

Each $M \in \Gr(kQ)$ has a largest submodule belonging to $\Fdim(kQ)$, namely
$$
\tau M:= \hbox{the  sum of all finite-dimensional graded submodules of $M$}.
$$
 
  \subsection{}
 Up to isomorphism and degree shift, the indecomposable projective graded left $kQ$-modules are  
 $$
 P_i=(kQ)e_i \cong kQ \otimes_{kI} ke_i, \qquad i \in I,
 $$
where $ke_i$ is the simple left $kI$-module 
 at vertex $i$. It follows that every projective module in $\Gr(kQ)$ is isomorphic 
 to $kQ \otimes_{kI} V$ for a suitable graded $kI$-module $V$.

 \begin{lem}
\label{lem.splitting}
Let $P,P' \in \Gr(kQ)$ be graded projective  modules generated by their degree $n$ components.
Every injective degree-preserving homomorphism  $f:P \to P'$ splits. 
\end{lem}
\begin{pf}
Without loss of generality we can assume $n=0$, $P=kQ \otimes_{kI} U$, and $P'=kQ \otimes_{kI} V$.
The natural map
$$
kQ \otimes_{kI} - : \Hom_{kI}(U,V) \to \Hom_{\Gr(kQ)}(kQ \otimes_{kI} U,kQ \otimes_{kI} V)
$$
is an isomorphism with inverse given by restricting a $kQ$-module homomorphism to 
the degree-zero components.
An injective  homomorphism $f:kQ \otimes_{kI} U \to kQ \otimes_{kI} V$ in $\Gr(kQ)$
 restricts to an injective $kI$-module homomorphism
$U \to V$  which splits because $kI$ is a semisimple ring.
\end{pf}

Part (1) of the next result is implied by \cite[Thm. 3.14]{AB} but because we only prove it for graded modules 
a simpler proof is possible.  
   
 \begin{prop}
 \label{prop.gr.coh}
$\phantom{xxx}$
 \begin{enumerate}
  \item 
Let $M$ be a finitely generated graded left $kQ$-module. Then $M$ is finitely presented if and only if  for all $n \gg 0$
 $$
 M_{\ge n} \cong  \bigoplus_{i \in I} P_i(-n)^{\oplus m_i}
 $$
 for some integers $m_i$ depending on $M$ and $n$. 
  \item 
  If $0 \to L \to M \to N \to 0$ is an exact sequence in $\gr(kQ)$, then 
   $0 \to L_{\ge n} \to M_{\ge n} \to N_{\ge n} \to 0$ splits for $n \gg 0$.
   \item{}
   If $\cM \in \qgr(kQ)$,  there is a projective $M \in \gr(kQ)$ such that $\cM \cong \pi^* M$.
   \item
    Every short exact sequence in $\qgr(kQ)$ splits. 
    \item
All objects in $\qgr(kQ)$ are injective and projective.  
\end{enumerate}
 \end{prop}
 \begin{pf}
 (1)
 ($\Leftarrow$) 
 Each $P_i$ is finitely presented because it is a finitely generated left ideal of the coherent ring $kQ$. 
Hence every $P_i(-n)$ is in $\gr(kQ)$. Therefore, if there is an integer $n$ such that $M_{\ge n}$ is a 
finite direct sum of various $P_i(-n)$s, then $M_{\ge n}$ is in $\gr(kQ)$ too. The hypothesis that $M$ is finitely generated implies that $M/M_{\ge n}$ is finite dimensional. But every finite dimensional graded $kQ$-module is 
a quotient of a direct sum of twists of the finitely presented finite dimensional module 
$kQ/kQ_{\ge 1}$ and is therefore in $\gr(kQ)$. 
In particular, $M/M_{\ge n} \in \gr(kQ)$. Since $\gr(kQ)$ is closed under extensions, $M$ is in $\gr(kQ)$ too.

 ($\Rightarrow$) 
Let $M \in \gr(kQ)$. Then there is an exact sequence $0 \to F' \stackrel{f}{\longrightarrow} F \to M \to 0$ in $\gr(kQ)$ with $F$ and $F'$ finitely generated graded projective $kQ$-modules. 
 Since $F'$, $F$, and $M$, are finitely generated, for all sufficiently large $n$ the modules
 $F'_{\ge n}$, $F_{\ge n}$, and $M_{\ge n}$, are generated as $kQ$-modules by $F'_n$, $F_n$, and $M_n$,
 respectively.  But $kQ$ is hereditary so $F_{\ge n}$ and $F'_{\ge n}$ are graded projective.
Now Lemma \ref{lem.splitting} implies that the restriction $f:F'_{\ge n}
 \to F_{\ge n}$ splits. Hence $M_{\ge n}$ is a direct summand of  $F_{\ge n}$. The result follows.

 (2) 
 By (1), $N_{\ge n}$ is projective for $n \gg 0$, hence the splitting.
 
 (3)
 There is some $M$ in $\gr(kQ)$ such that $\cM \cong \pi^*M$. But $\pi^* M \cong \pi^*(M_{\ge n})$ for all $n$
 so (3) follows from (1). 
 
 (4)
 By \cite[Cor. 1, p. 368]{Gab}, every short exact sequence in $\qgr(kQ)$  is of the form
\begin{equation}
\label{eq.qgr.ses}
0 \longrightarrow  \pi^* L \stackrel{\pi^*f}{\longrightarrow}  \pi^* M \stackrel{\pi^*g}{\longrightarrow}  \pi^* N \longrightarrow 0
\end{equation}
for some exact sequence $0 \longrightarrow  L \stackrel{f}{\longrightarrow}  M \stackrel{g}{\longrightarrow}   N \longrightarrow 0$  in $\gr(kQ)$. But (\ref{eq.qgr.ses}) is also obtained by applying $\pi^*$ to the restriction
 $0 \to L_{\ge n} \to M_{\ge n} \to N_{\ge n} \to 0$ which splits for $n \gg 0$. Hence  (\ref{eq.qgr.ses}) splits.
  \end{pf}
  
\subsection{}
By Proposition \ref{prop.gr.coh}(5), $\cO$ is a projective object in $\qgr(kQ)$.
 
\begin{lemma}
$\cO$ is a projective object in $\QGr(kQ)$.
\end{lemma}
\begin{pf}
As noted in the proof of Proposition \ref{prop.gr.coh}(4), an epimorphism in $\QGr(kQ)$ is 
necessarily of the form $\pi^*g:\pi^*M \to \pi^*N$ for some surjective homomorphism $g:M \to N$ in 
$\Gr(kQ)$. 
Let $\eta:\cO \to \pi^*N$ be a morphism in $\QGr(kQ)$. Then
$$
\eta \in \liminj \Hom_{\Gr(kQ)}(kQ_{\ge n}, N/N')
$$
where the direct limit is taken over all $n \in \NN$ and all $N' \subset N$ such that $N/N'$ is 
the sum of its finite dimensional submodules, so $\eta=\pi^* h$ for some $n$, some $N'$, and some $h:kQ_{\ge n} \to N/N'$. 
Since $kQ_{\ge n}$ is a projective object in $\Gr(kQ)$, 
$h$ factors through $N$ and for the same reason $h$ factors through $g$. Hence there is a morphism $\gamma:\cO \to \pi^* M$ such that $\eta=(\pi^*g) \circ \gamma$.  
\end{pf}

\begin{lemma}
$\Hom_{\QGr(kQ)}(\cO,-)$ commutes with all direct sums in $\QGr(kQ)$.
\end{lemma}
\begin{pf}
Let 
$$
\cM= \bigoplus_{\lambda  \in \Lambda} \cM_\lambda 
$$
be a direct sum in $\QGr(kQ)$.  Let $M_\l$, $\l \in \L$, be graded $kQ$-modules
such that $\pi^*M_\lambda =\cM_\l$. Because $\pi^*$ has a right adjoint it commutes with direct sums. 
Hence $\cM=\pi^* M$
where $M=\oplus_{\lambda  \in \Lambda} M_\lambda $. 
Because $kQ_{\ge n}$ is a finitely generated module we obtain the second equality in the computation
\begin{align*}
\Hom_{\QGr(kQ)}(\cO,\cM) & = \liminj_n \Hom_{\Gr(kQ)} (kQ_{\ge n}, \oplus_{\lambda  \in \Lambda} M_\l) \\
&
=   \liminj_n \bigoplus_{\lambda  \in \Lambda} \Hom_{\Gr(kQ)} (kQ_{\ge n},M_\lambda ) 
\\
& = \bigoplus_{\lambda  \in \Lambda}   \liminj_n \Hom_{\Gr(kQ)} (kQ_{\ge n},M_\lambda ) 
\\
 & = \bigoplus_{\lambda  \in \Lambda}  \Hom_{\QGr(kQ)} (\cO,\cM_\lambda ).
\end{align*}
This proves the lemma. 
\end{pf}

\subsubsection{Notation}
We write $\cP_i=\pi^*P_i$ for the images of the indecomposable projectives in $\QGr(kQ)$.

\subsubsection{}
If $\cS$ is a set of objects in an additive category $\sA$ we write $\add(\cS)$ for the smallest 
full subcategory of $\sA$ that contains $\cS$ and is closed under direct summands and finite direct sums.

\begin{lem}
\label{lem.O.genor}
For every positive integer $n$, $\cO(-n)\in \add(\cO)$.
\end{lem}
\begin{pf} 
Let
\begin{align*}
I_m:= & \{ v \in I \; | \; \hbox{there is a path of length $m$ that ends at $v$}\}, \; m \ge 1,
\\
I_0:= & I-I_1,
\\
I_\infty:= &\bigcap_{m=1}^\infty I_m,
\\
\sT_m:= & \add\{\cP_i \; | \; i \in I_m\}, \; 0 \le m \le \infty.
\end{align*}
The vertices in $I_0$ are the sources.
A vertex is in $I_\infty$ if and only if for every $m \ge 1$ there is a path of length $m$ ending at it.
If $m \gg 0$, then
$$
I_1 \supset I_2 \supset \cdots \supset I_m=I_{m+1}=\cdots = I_\infty
$$
and, consequently, 
$$
\sT_1 \supset \sT_2 \supset \cdots \supset \sT_m=\sT_{m+1}=\cdots = \sT_\infty.
$$
 
 To prove the lemma it suffices to show that $\cO(-1)\in \add(\cO)$ because, if it is, an induction argument would complete the proof: $\cO(-1)\in \add(\cO)$ implies that $\cO(-2)\in \add(\cO(-1)) \subset \add(\cO)$,
 and so on. But $\cO(-1) = \oplus_{i\in I} \cP_i(-1)$ is a direct sum of an object in $\sT_0(-1)$ and an object in 
 $\sT_1(-1)$ so it suffices to show that $\add(\cO)$ contains $\sT_0(-1)$ and $\sT_1(-1)$.  

If $j$ is a sink, then $P_j=ke_j$ so $\cP_j=0$. 

Suppose $j \in I_m$ and $j$ is not a sink.
There is an exact sequence
 $$
 0 \to \bigoplus_{a \in s^{-1}(j)} P_{t(a)}(-1) \stackrel{(\cdot a)}{\longrightarrow} P_j \to ke_j \to 0
 $$
 where $ke_j$ denotes the simple module concentrated at vertex $j$. Therefore
 $$
 \cP_j \cong \bigoplus_{a \in s^{-1}(j)} \cP_{t(a)}(-1)
 $$
 for every vertex $j$. If $a \in s^{-1}(j)$, then $t(a) \in I_{m+1}$. Therefore $\sT_m \subset \sT_{m+1}(-1)$. On the other hand, if $m \ge 1$ and $i \in I_{m+1}$, there is an arrow $a$ such that $t(a)=i$ and $s(a) \in I_m$ so
  $\cP_i(-1)$ is a direct summand of $\cP_{s(a)}$. Hence $\sT_{m+1}(-1) \subset \sT_m$. 
 
The previous paragraph shows that $ \sT_0 \subset \sT_1(-1)$ and
$
 \sT_m = \sT_{m+1}(-1)
$
 for all $m \ge 1$.
Thus
 \begin{equation}
 \label{eq.twist.T}
 \sT_1=\sT_2(-1) = \cdots = \sT_m(-m+1)
 \end{equation}
 for all $m \ge 1$. For $m \gg 0$, $\sT_m=\sT_{m+1}$ so
 $
 \sT_m(-1)=\sT_{m+1}(-1)=\sT_m$.
Hence, for $m \gg 0$,  $\sT_m=\sT_m(n)$ for all $n \in \ZZ$. Therefore (\ref{eq.twist.T}) implies $\sT_1=\sT_1(n)$ for all $n \in \ZZ$.

Since $\cO = \oplus_{i\in I} \cP_i$, 
$\sT_1\subset \add(\cO)$. Therefore $\sT_1(-1) =\sT_1 \subset \add(\cO)$ and 
$\sT_0(-1) \subset \sT_1(-2) =\sT_1(-1) \subset \add(\cO)$. 
\end{pf}

\begin{prop}
\label{prop.O.genor}
$\qgr(kQ)=\add(\cO)$.
\end{prop}
\begin{pf}
Let $\cM \in \qgr(kQ)$. There is some $M$ in $\gr(kQ)$ such that $\cM \cong \pi^*M$. But $\pi^* M \cong \pi^*(M_{\ge n})$ for all $n$ so, by Proposition \ref{prop.gr.coh}(1), if $n \gg 0$ there are  integers $m_i$ such that
$$
\cM \cong \bigoplus_{i \in I} \cP_i(-n)^{\oplus m_i}.
$$
Each $\cP_i(-n)$ belongs to $\add(\cO)$ by Lemma \ref{lem.O.genor} so $\cM \in \add(\cO)$. 
\end{pf}

  \begin{thm}
  Let $\Mod S(Q)$ be the category of right $S(Q)$-modules.
 The functor $\Hom_{\QGr(kQ)}(\cO,-)$ provides an equivalence of categories
  $$
\QGr(kQ) \equiv \Mod S(Q)
 $$
 that sends $\cO$ to $S(Q)$.  This equivalence restricts to an equivalence between $\qgr(kQ)$ and $\mod S(Q)$, 
the category of finitely presented $S(Q)$-modules. 
 \end{thm}
 \begin{pf}
By Proposition \ref{prop.O.genor}, $\cO$ is a generator in $\qgr(kQ)$. Every object in $\QGr(kQ)$
is a direct limit of objects in $\qgr(kQ)$ so $\cO$ is also a generator in $\QGr(kQ)$.   
Since $\cO$ is a finitely generated, projective generator in  the Grothendieck category $\QGr(kQ)$,
$$
\Hom_{\QGr(kQ)}(\cO,-): \QGr(kQ) \to \Mod \big(\End_{\QGr(kQ)}(\cO)\big).
$$
is an equivalence of categories \cite[Example X.4.2]{Sten}.  
The result now follows from the isomorphism 
 $\End_{\QGr(kQ)}( \cO) \cong S(Q)$ in Theorem \ref{thm.End.O}.
 \end{pf}


\section{Sinks and sources can be deleted}
\label{sect.sink.source}

\subsection{}
A vertex is a {\sf sink} if no arrows begin at it and a 
{\sf source} if no arrows  end at it.

\begin{thm}
\label{thm.source.sink}
 Suppose quivers $Q$ and $Q'$ become the same after repeatedly removing
sources and sinks and attached arrows. Let $(Q^\circ,I^\circ)$ be the quiver without sources or sinks that is 
obtained from $Q$ by this process. Then there is an equivalence of categories
$$
\QGr(kQ) \equiv \QGr(kQ^\circ) \equiv \QGr(kQ').
$$
The equivalence of categories is induced by sending a representation \newline  $(M_i,M_a; i \in I, a \in Q_1)$
of $Q$ to the representation $(M_i,M_a; i \in I^\circ, a \in Q^\circ_1)$ of $Q^\circ$.  
A quasi-inverse to this is induced by the functor that sends 
a representation $(M_i,M_a; i \in I^\circ, a \in Q^\circ_1)$ of $Q^e$ to the representation $(M_i,M_a; i \in I, a \in Q_1)$
of $Q$ where $M_i=0$ if $i \notin I^\circ$ and $M_a =0$ if $a \notin Q^\circ_1$. 
\end{thm}

\subsection{}
\label{sect.no.sss}
The fact that sinks and sources can be deleted is reminiscent of three other results in the literature. 

The category $\QGr(kQ)$ is related to the dynamical system with topological space 
the bi-infinite paths in $Q$ viewed as a subspace of $Q_1^\ZZ$ and automorphism the edge shift $\s$ defined by $\s(f)(n)=f(n+1)$. No arrow that begins at a source and no arrow that ends at a sink appears in any bi-infinite path so, as  remarked after Example 2.2.8 in \cite{LM}, since $Q^\circ$ ``contains the only part of $Q$ used for symbolic dynamics, we will usually confine our attention to [quivers such that $Q=Q^\circ$]''. 
See also \cite[Prop. 2.2.10]{LM}. 

Second, as remarked on page 18 of \cite{R}, ``Cuntz-Krieger algebras are the C$^*$-algebras of finite graphs 
with no sinks or sources''. 

Third, in section 4 of \cite{XWChen} it is shown that the singularity category of an artin algebra is not changed by 
deleting or adding sources or sinks.

\subsection{}
Theorem \ref{thm.source.sink} follows from the next two results.

\begin{prop}
\label{prop.sink}
Let $t$ be a sink in $(Q,I)$ and  $Q'$ the quiver with vertex set $I':=I-\{t\}$ and arrows 
$$
Q'_1:=\{\hbox{arrows in $Q$ that do not end at $t$}\}.
$$
Then the functor $i_*:\Gr(kQ') \to \Gr(kQ)$ that sends a representation of $Q'$ to the ``same'' representation of 
$kQ$ obtained by putting 0 at vertex $t$ induces an equivalence of categories
$$
\QGr(kQ') \equiv \QGr(kQ)
$$
that sends $\cO'$ to $\cO$.
\end{prop}
\begin{pf}
Since $t$ is a sink, $(kQ)e_t=ke_t$. Hence  $e_t(kQ)$ is a two-sided ideal of $kQ$ and, if $M$ is a left $kQ$ module, then $e_tM$ is a submodule of $M$. 

It is clear  that $i_*$ is the forgetful functor induced by the
homomorphism $kQ \to kQ/(e_t) = kQ'$. The functor $i_*$ has a left adjoint $i^*$ and a right adjoint $i^!$.
The functor $i^*$ sends a $kQ$-module $M$ to $M/e_tM$.
It is easy to see that the counit $i^*i_* \to  \id_{\Gr(kQ')} $ is an isomorphism.

Both $i_*$ and $i^*$ are exact.

Since $i_*$ and $i^*$ send direct limits of finite dimensional modules to
direct limits of  finite dimensional modules they induce
functors $\iota_*:\QGr(kQ') \to \QGr(kQ)$ and $\iota^*:\QGr(kQ) \to \QGr(kQ')$.  Because 
$i^*i_* \cong \id_{\Gr(kQ')} $ we have $\iota^*\iota_* \cong \id_{\QGr(kQ')}$.  

If $M \in \Gr(kQ)$ there is an exact sequence $0 \to e_t M \to M \to i_*i^*M \to 0$. 
If $a$ is any arrow, then $ae_t=0$. Therefore $e_tM$ is a direct sum of 1-dimensional left $kQ$-modules,
hence in $\Fdim(kQ)$. It follows that the unit $\id_{\Gr(kQ)} \to i_*i^*$
induces an isomorphism $\id_{\QGr(kQ)} \cong \iota_*\iota^*$.  

Hence $\iota_*$ and $\iota^*$ are mutually quasi-inverse equivalences.

Since $kQ'=kQ/e_tkQ$, $i_*$ sends $kQ'$ to $kQ/e_tkQ$. The natural homomorphism 
$kQ \to kQ/e_tkQ=i_*(kQ')$ becomes an isomorphism in $\QGr(kQ)$ 
because $e_tkQ \in \Fdim(kQ)$.  Hence $\iota_*\cO'=\cO$.
\end{pf}

\begin{prop}
\label{prop.source}
Let $s$ be a source in $(Q,I)$ and  $Q'$ be the quiver with vertex set $I':=I-\{s\}$ and arrows 
$$
Q'_1:=\{\hbox{arrows in $Q$ that do not begin at $s$}\}.
$$
The functor $i_*:\Gr(kQ') \to \Gr(kQ)$ that sends a representation of $Q'$ to the ``same'' representation of $kQ$  obtained by putting 0 at vertex $s$ induces an equivalence, $\iota_*$, of categories
$$
\QGr(kQ') \equiv \QGr(kQ).
$$
Furthermore, $\cO \cong \iota_*\cO' \oplus \cP_s$ where $\cP_s=\pi^*(kQe_s)$ and $\pi^*$
is the quotient functor $\Gr(kQ) \to \QGr(kQ)$.
\end{prop}
\begin{pf}
Every $kQ'$-module becomes a $kQ$-module through the homomorphism $\varphi:kQ \to kQ/(e_s) = kQ'$; 
this is the exact fully faithful embedding $i_*$. A right 
adjoint to $i_*$ is given by the functor $i^!$,
$$
i^!M:=\Hom_{kQ}(kQ',M) = \{m \in M \; | \; e_sm=0\}=(1-e_s)M.
$$
It is clear that the unit $\id_{\Gr(kQ')}  \to i^!i_*$ is an isomorphism of functors.

Both $i_*$ and $i^!$ are exact.

Since $i_*$ and $i^!$ send direct limits of finite dimensional modules to
direct limits of  finite dimensional modules there  are unique functors 
$\iota_*:\QGr(kQ') \to \QGr(kQ)$ and $\iota^!:\QGr(kQ) \to \QGr(kQ')$ such that the diagrams
$$
\UseComputerModernTips
\xymatrix{
\Gr(kQ') \ar[r]^{i_*} \ar[d] & \Gr(kQ)  \ar[d] 
\\
\QGr(kQ') \ar[r]_{\iota_*} & \QGr(kQ) 
}
\qquad \hbox{and} \qquad 
\UseComputerModernTips
\xymatrix{
\Gr(kQ) \ar[r]^{i^!} \ar[d] & \Gr(kQ') \ar[d] 
\\
\QGr(kQ) \ar[r]_{\iota^!} & \QGr(kQ')
}
$$
commute \cite[Sect. III.1]{Gab}. (The vertical arrows in the diagrams are the quotient functors.)

Because  $\id_{\Gr(kQ')} \cong i^!i_*$, $\iota^!\iota_* \cong \id_{\QGr(kQ')}$.  


If $M \in \Gr(kQ)$, there is an exact sequence $0 \to i_*i^!M \to M \to \Mol \to 0$ in which $\Mol$ is supported only at the vertex $s$; a module supported only  at $s$ is a sum of 1-dimensional $kQ$-modules
so belongs to $\Fdim(kQ)$.  It follows that the unit $i_*i^! \to \id_{\Gr(kQ)}$
induces an isomorphism $\id_{\QGr(kQ)} \cong \iota_*\iota^!$.  

Hence $\iota_*$ and $\iota^!$ are mutually quasi-inverse equivalences.

The isomorphism $\cO \cong \iota_*\cO' \oplus \cP_s$ is proved in section \ref{sect.end.pf} below.
\end{pf}

\subsection{}
It need not be the case that the equivalence $\iota_*$ in Proposition \ref{prop.source} sends $\cO'$ to $\cO$. 
We will show that $\iota_*\cO' \not\cong \cO$ for the quivers
$$
Q=  \qquad  \UseComputerModernTips
\xymatrix{
 s  \ar[r]^a & v   \ar@(ur,dr)[]^b
 }
 \qquad \hbox{and} \qquad 
 Q'= \qquad  \UseComputerModernTips
\xymatrix{
 v \ar@(ur,dr)[]^b  
 }
$$
Since $kQ'$ is a polynomial ring in one variable, $\QGr(kQ') \equiv \Mod k$ and this equivalence sends $\cO'$, the image of $kQ'$ in $\QGr(kQ')$, to $k$. Since $\iota_*$ is an equivalence it follows that $\iota_* \cO'$ is indecomposable. Now $kQ$ is isomorphic as a graded left $kQ$-module to the direct sum of the projectives $P_s$ and $P_v$. Right multiplication by the arrow $a$ gives an isomorphism $P_v \to (P_s)_{\ge 1}(1)$
in $\Gr(kQ)$. Let $\cP_s=\pi^* P_s$ and $\cP_v=\pi^* P_v$. Then
$$
\cO = \cP_v \oplus \cP_s \cong \cP_v \oplus \cP_v(-1).
$$
Hence $\cO \not\cong \iota_*\cO'$. 

Right multiplication by $b$ induces an isomorphism $\cP_v \stackrel{\sim}{\longrightarrow} \cP_v(-1)$. 

\subsection{}
\label{sect.end.pf}
We now prove the last sentence of Proposition \ref{prop.source}. 

Because $s$ is a source, the two-sided ideal $(e_s)$ is equal to $kQe_s$.
Hence as a graded left $kQ$-module $i_*(kQ')$ is isomorphic to $kQ/kQe_s$ which is isomorphic to 
$kQ(1-e_s)$. The claim that $\cO \cong \iota_*\cO' \oplus \cP_s$ now follows from the decomposition $kQ=kQ(1-e_s) \oplus kQe_s$.

\subsubsection{}
Because $\cO \cong \iota_*\cO' \oplus \cP_s$, $\iota^!\cO \cong \cO' \oplus \iota^!\cP_s$. 
Moreover, $\iota^!\cP_s$ is isomorphic to $\oplus_{s(a)=s} \cP'_{t(a)}(-1)$ where $\cP'_{t(a)}$
is the image in $\QGr(kQ')$ of $(kQ')e_{t(a)}$. 


\section{Description of $S(Q)$}
\label{sect.brat}

We will give two different descriptions of $S(Q)$. 

In section \ref{sect.B.diag}, we describe $S(Q)$ in terms of its Bratteli diagram. See \cite{B} and \cite{Eff} for
information about Bratteli diagrams.

In section \ref{sect.leavitt}, we show that when $Q$ has no sinks or sources $S(Q)$ is isomorphic to the degree zero component of the Leavitt path algebra, $L(Q)$, associated to $Q$
and, because $L(Q)$ is strongly graded,
$$
\QGr(kQ) \equiv \Gr_{\ell} L(Q) \equiv \Mod_{\ell} L(Q)_0
$$
where the subscript $\ell$ means {\it left} modules.

It is well-known that $L(Q)$ is a dense subalgebra of the Cuntz-Krieger algebra ${\mathscr O}_Q$,  associated to $Q$. The philosophy of non-commutative geometry suggests that $kQ$ is a homogeneous coordinate ring for a non-commutative scheme whose underlying non-commutative topological space has ${\mathscr O}_Q$ as its ring of ``continuous $\CC$-valued functions''.

\subsection{The Bratteli diagram for $S=S(Q)$}
\label{sect.B.diag}
Because $kI$ has $|I|$ isoclasses of simple modules and $S_n$ is the endomorphism 
ring of a left $kI$-module, $S_n$ is a product of at most $|I|$ matrix algebras of various sizes;
although fewer than $|I|$ matrix algebras may occur in the product it is better to think  there are $|I|$ of 
them with the proviso that some (those corresponding to sources) might be $0 \times 0$ matrices.

The $n^{\th}$ level of the Bratteli diagram for $S(Q)$ therefore consists of $|I|$ vertices, each denoted by 
$\bullet$, that we label $(n,i)$,
$i \in I$. The vertex labelled $(n,i)$ represents the summand $\End_{kI}(e_i(kQ))$ of $S_n$; this endomorphism 
ring is isomorphic to a matrix algebra $M_r(k)$ for some integer $r$; it is common practice to replace the symbol $\bullet$ at $(n,i)$ by the integer $r$; we then say that $r$ is {\sf the number at vertex } $(n,i)$ . We do this for some of the examples in section \ref{sect.egs}. 

We will see that the number of edges from $(n,i)$ to $(n+1,j)$ is the same as the number of edges 
from $i$ to $j$ in $Q$. The Bratteli diagram is therefore stationary in the terminology of \cite{E}, and the following example, which appears in \cite{HS}, illustrates how to pass from the quiver to the associated 
Bratteli diagram.
$$
   \UseComputerModernTips
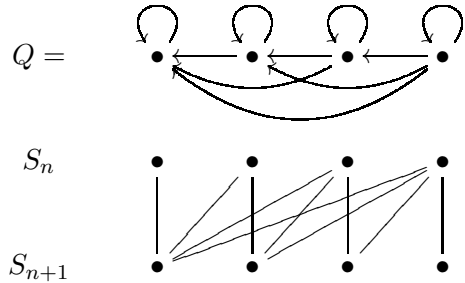
\xymatrix{
\\
Q= & \bullet   \ar@(ur,ul)[]  & \ar[l] \bullet   \ar@(ur,ul)[]  & \ar[l]\bullet   \ar@(ur,ul)[]  \ar@/^1pc/[ll]
 & \ar[l]\bullet   \ar@(ur,ul)[]  \ar@/^1pc/[ll] \ar@/^2pc/[lll]
 \\
S_n &  \bullet   \ar@{-}[d] & \ar@{-}[dl] \bullet   \ar@{-}[d]  & \ar@{-}[dl]\bullet   \ar@{-}[d]  \ar@{-}[dll]
 & \ar@{-}[dl]\bullet   \ar@{-}[d]  \ar@{-}[dll] \ar@{-}[dlll]
 \\
S_{n+1} &  \bullet  &\bullet   &  \bullet   &  \bullet 
}
$$

 \subsection{}
Let $C := (c_{ij})_{i,j \in I}$ be the incidence matrix for $Q$ with the convention that
$$
c_{ij}=\hbox{the number of arrows from $j$ to $i$}.
$$
The $ij$-entry in $C^n$, which we denote by $c^{[n]}_{ij}$, 
is the number of paths of length $n$ from $j$ to $i$.
The number of paths of length $n$ ending at vertex $i$ is  
$$
p_{n,i}:= \sum_{j \in I} c^{[n]}_{ij}.
$$

\begin{prop}
\label{prop.brat} 
The sum of the left $kI$-submodules of $kQ_n$ isomorphic to $E_i$ is equal to $e_i(kQ_n)$.
Its dimension is equal to $p_{n,i}$ and 
\begin{equation}
\label{eq.Sn.isom}
S_n \cong \bigoplus_{i \in I} \End_{kI} \big(e_i(kQ_n)\big)    \cong \bigoplus_{i \in I} M_{p_{n,i}}(k).
\end{equation}
Referring to the Bratteli diagram for $\liminj S_n$, 
the number  at the vertex labelled $(n,i)$ is $p_{n,i}$, and the 
number of edges from $(n,i)$ to $(n+1,j)$ is $c_{ji}$. 

Composing with the inclusions and projections in (\ref{eq.Sn.isom}), 
the components of the map $\theta_n:S_n \to S_{n+1}$ 
are the maps
$$
 \End_{kI} \big(e_i(kQ_n)\big) \to S_n \stackrel{\theta_n}{\longrightarrow} S_{n+1} \to 
  \End_{kI} \big(e_j(kQ_{n+1})\big) 
$$
that send a matrix in $\End_{kI} \big(e_i(kQ_n)\big)$ to $c_{ji}$ ``block-diagonal'' copies of itself in $ \End_{kI} \big(e_j(kQ_{n+1})\big)$.
\end{prop}
\begin{pf}
The irreducible representation of $Q$ at vertex $i$ is $E_i = ke_i$.  
Since a path $p$ ends at vertex $i$ if and only if $p=e_ip$,
the multiplicity of $E_i$ in a composition series for $kQ_n$ as a left $kI$-module   is 
\begin{align*}
[kQ_n:E_i]  &= \dim_ke_i(kQ_n) 
\\
& = \hbox{the number of paths of length $n$ ending at $i$}
\\
& =p_{n,i}.
\end{align*}
The existence of the left-most isomorphism in (\ref{eq.Sn.isom}) follows at once from the fact that $kI$ is a semisimple ring; the second isomorphism follows from the analysis in the first part of this paragraph.

Up to isomorphism, $\{E_i^*=\Hom_k(E_i,k) \; | \; i \in I\}$ is a complete set of simple right $kQ$ modules.
It follows that the $kI$-bimodules 
$$
E_{ij}:=E_i \otimes E_j^*, \qquad i,j \in I,
$$ 
form a complete set of isoclasses of simple $kI$-bimodules. If $a$ is an arrow from $j$ to $i$
there is a $kI$-bimodule isomorphism $ka \cong E_{ij}$ so the multiplicity of $E_{ij}$ in $kQ_1$ is the number of arrows from $j$ to $i$, i.e.,   $[kQ_1: E_{ij}] =  c_{ij}$. More explicitly, $E_{ij}^{\oplus c_{ij}} \cong e_i(kQ_1)e_j$. 

The image in $S_{n+1}$ of a map $f \in S_n$  is the map $kQ_1 \otimes f$. Hence if $f$ is belongs to the component  $\End_{kI} \big(e_i(kQ_n)\big)$ of $S_n$, the component of $kQ_1 \otimes f$ in 
$  \End_{kI} \big(e_j(kQ_{n+1})\big)$
is $e_jkQ_1e_i \otimes f$. But the dimension of $e_jkQ_1e_i$ is $[kQ_1: E_{ji}]=c_{ji}$. 
\end{pf}

\subsection{Leavitt path algebras}
\label{sect.leavitt}
Goodearl's survey \cite{KG} is an excellent introduction to Leavitt path algebras. 

Since $\QGr(kQ)$ is unchanged when $Q$ is replaced by the quiver obtained by repeatedly deleting 
sources and sinks,  the essential case is when $Q$ has no sinks or sources. 

For the remainder of section \ref{sect.brat} we assume $Q$ has no sinks or sources.
This is equivalent to the hypothesis that $Q=Q^\circ$.

\subsubsection{}
\label{ssect.1-3}
Under the hypothesis that $Q=Q^\circ$,
\begin{enumerate}
 \item 
 the Leavitt path algebra of $Q$,  $L(Q)$, is a universal localization of 
 $kQ$ in the sense of \cite[Sect. 7.2]{Cohn-FRR} or \cite[Ch. 4]{Scho};
  \item 
 $L(Q)$ is a strongly $\ZZ$-graded ring and $L(Q)_0 \cong S(Q)^{\op}$ so 
   $$
   \Gr_\ell  L(Q) \equiv \Mod_\ell L(Q)_0 \equiv \Mod_r S(Q) \equiv \QGr(kQ)
   $$
   where the subscripts $\ell$ and $r$ denote left and right modules.
   \item{}
$L(Q)$ is a dense subalgebra of the Cuntz-Krieger algebra ${\mathscr O}_Q$.
\end{enumerate}

\subsubsection{}
Statement (1)  in \ref{ssect.1-3} holds for {\it every} finite $Q$ (i.e., the hypothesis $Q=Q^\circ$ is not needed) \cite{AB}.  For \ref{ssect.1-3}(3), see \cite{AT} and \cite{R}.

The fact that $L(Q)$ is strongly graded when $Q=Q^\circ$ is proved in  \cite[Thm. 3.11]{H}. We give an alternative proof of this in Proposition \ref{prop.LQ+-}(6). 
A $\ZZ$-graded ring $R$ is {\sf strongly graded} if $R_nR_{-n}=R_0$ for all $n$. 
When $R$ is strongly graded $\Gr R$ is equivalent to $\Mod R_0$ via the functor $M \rightsquigarrow M_0$ \cite[Thm. 2.8]{Dade}.  That explains the left-most equivalence in (2).  
The second equivalence in \ref{ssect.1-3}(2) follows from the fact that $S(Q)^{\op} \cong L(Q)_0$ which we will prove in 
Theorem \ref{thm.L+S} below.

\subsection{}
Let $E_i$ be the 1-dimensional left $kQ$-module supported at vertex $i$ and concentrated in degree zero.
Because $Q$ has no sinks $E_i$ is not projective and its minimal projective resolution is
\begin{equation}
\label{eq.min.res.Si}
0 \longrightarrow \bigoplus_{a \in s^{-1}(i)} P_{t(a)} \stackrel{f_i}{\longrightarrow} P_i \longrightarrow E_i \longrightarrow 0
\end{equation}
where $P_j=(kQ)e_j$ and the direct sum is over all arrows starting at $i$. 
Elements in the direct sum  will be written as row vectors $(x_a,x_b,\ldots)$ with 
 $x_a \in P_{t(a)}$, $x_b \in P_{t(b)}$, and so on. The map $f_i$ 
is right multiplication by the column vector $(a,b,\ldots)^\sT$ where $a,b,\ldots$ are the arrows starting at $i$, i.e.,
\begin{equation}
\label{eq.fi}
f_i(x_a,x_b,\ldots ) = (x_a,x_b,\ldots ) \begin{pmatrix} a \\ b \\ \vdots\end{pmatrix} =x_aa+x_bb+\cdots \, \in \, (kQ)e_i.
\end{equation}

\subsubsection{Definition of $L(Q)$ as a universal localization}
We refer the reader to \cite[Sect. 7.2]{Cohn-FRR} and \cite[Ch. 4]{Scho} for details about universal localization.

Let $\Sigma = \{f_i \; | \; i \in I\}$ and let $$L(Q):=\Sigma^{-1}(kQ)$$ be the universal localization of $kQ$ at $\Sigma$. Since $Q$ will not change in this section we will often write $L$ for $L(Q)$.

If $x \in kQ$ we continue to write $x$ for its image in $L$ under the universal $\Sigma$-inverting map 
$kQ \to L$. Let
$$
P':= \bigoplus_{a \in s^{-1}(i)} Le_{t(a)}.
$$
The defining property of $L$ is that the map $kQ \to L$ is universal subject to the
 condition that  applying $L \otimes_{kQ} - $ to  (\ref{eq.min.res.Si}) produces an isomorphism
$$
   \UseComputerModernTips
\xymatrix{
\id_L \otimes f_i :P'  \ar[rr]^>>>>>>>>>>>{\hdot \begin{pmatrix} a \\ b \\ \vdots\end{pmatrix} } &&  Le_i
}
$$
for all $i$. Thus $L \otimes_{kQ} E_i=0$ for all $i \in I$.

Every $L$-module homomorphism $Le_i \to Le_j$ is right multiplication by an element of $L$ so
 \cite[Thm. 4.1]{Scho} tells us $L$ is generated by $kQ$ and elements $a^*,b^*,\ldots$ such that
the inverse of $\id_L \otimes f_i$ is right multiplication by the row vector $(a^*,b^*,\ldots)$ where
 $a^*= e_ia^*e_{t(a)}$, etc. 
In particular, the defining relations for $L$ are given by
 $$
 (a^*,b^*,\ldots) \begin{pmatrix} a \\ b \\ \vdots\end{pmatrix} = \id_{Le_i} 
  \qquad \hbox{and} \qquad
 \begin{pmatrix} a \\ b \\ \vdots\end{pmatrix}  (a^*,b^*,\ldots) = \id_{P'}. 
 $$
Since $\id_{Le_i}$ is right multiplication by $e_i$ and 
$\id_{P'} $ is right multiplication by 
$$
\begin{pmatrix}
  e_{t(a)} & 0 & 0    & \cdots & 0    \\
    0   & e_{t(b)} & 0 & \cdots & 0   \\
    \vdots &&& & \vdots   
\end{pmatrix},
$$ 
$L=kQ \, \langle a^* \; | \; a \in Q_1\rangle$ modulo the relations  
\begin{align*}
e_{s(a)}a^*e_{t(a)} = & \, a^* \qquad \phantom{xxxtxix} \hbox{for all arrows $a \in Q_1$,}
\\
 aa^* = & \, e_{t(a)}  \qquad \phantom{xxxxx} \hbox{for all arrows $a \in Q_1$,}
\\
ab^* = & \,  0 \qquad  \phantom{xxxxxixi}  \hbox{if $a$ and $b$ are different arrows,}
\\
e_i= & \sum_{a \in s^{-1}(i)} a^*a  \qquad \hbox{for all $i \in I$} .
\end{align*}

\subsubsection{}
Our $L(Q)$ is not defined in the same way as the algebra $L_k(Q)$
defined in \cite[Sect. 1]{KG}.
Because our notational convention for composition of paths is the reverse of that in \cite[Sect. 1.1]{KG}
the relations for $L(Q)$ just above  are the opposite of those for $L_k(Q)$ in \cite[Sect. 1.4]{KG}.
(If we had defined $L(Q)$ by inverting homomorphisms between {\it right} instead of left modules we would have obtained the relations in \cite[Sect. 1.4]{KG} but then our convention for  composition of paths would have created the problems discussed at the end of \cite[Sect. 1.8]{KG}.)
As a consequence, our $L(Q)$ is anti-isomorphic to $L_k(Q)$. 
{\it However,} as explained at the end of \cite[Sect. 1.7]{KG},  $L_k(Q)$ is anti-isomorphic to itself via 
a map that sends each arrow $a$ to $a^*$ and fixes $e_i$ 
for each vertex $i$.  Thus, our $L(Q)$ {\it is} isomorphic to the algebra $L_k(Q)$ defined in \cite{KG}. 

\begin{prop}
The algebra $L(Q)=\Sigma^{-1}(kQ)$ is isomorphic and anti-isomorphic to the Leavitt path algebra $L_k(Q)$ defined in 
\cite[Sect. 1.4]{KG}.
\end{prop}

\subsubsection{}
Our convention for composition of paths is that used by (most of) the finite dimensional algebra community 
and by  Raeburn \cite[Rmk. 1.1.3]{R}. 
However, it is not the convention adopted in \cite{AT} (see \cite[Defn. 2.9]{AT}).

\subsection{}
If $a,b,\ldots,c,d \in Q_1$ and $p=dc\ldots ba$ we define $p^*:=a^*b^*\ldots c^* d^*$. 
If $p$ and $q$ are paths of the same length, then
\begin{equation}
\label{LQ.products}
pq^*=\d_{p,q} e_{t(q)} =\d_{p,q} e_{t(p)}.
\end{equation}

If $p$ is a path in $Q$ of length $n$ we write $|p|=n$. 
We give $L(Q)$ a $\ZZ$-grading by declaring that
$$
\deg a =1 \quad \hbox{and} \quad \deg a^* =-1 \quad \hbox{for all $a \in Q_1$}.
$$

For completeness we include the following well-known fact.

\begin{lem}
\label{lem.LQn}
\label{lem.LQ.span}
The degree-$n$ component of $L(Q)$ is 
$$
L_n={\rm span}\{p^*q \; | \; \hbox{$p$ and $q$ are paths such that $|q|-|p|=n$}\}.
$$
\end{lem}
\begin{pf}
Certainly  $L(Q)$ is spanned by words in the letters $a$ and $a^*$, $a \in Q_1$. Let $w$ be a non-zero word and $ab^*$ a subword of $w$  with $a,b \in Q_1$. Since $w \ne 0$, $ab^*=aa^*=e_{t(a)}$; but $e_{t(a)}$ can be absorbed into the letters on either side of $aa^*$ so, repeating this if necessary, $w=p^*q$ for some paths $p$ and $q$. The degree of $p^*q$ is $|q|-|p|$ so the result follows.
\end{pf}

\begin{thm}
\label{thm.L+S}
The algebras $L(Q)_0$ and $S(Q)$ are anti-isomorphic, 
$$
L(Q)_0 \cong S(Q)^{\op}.
$$
\end{thm}
\begin{pf}
By definition, $S(Q)$ is the ascending union of its subalgebras $S_n=\End_{kI}( kQ_n)$. 

We will sometimes write $L$ for $L(Q)$.

It is clear that $L_0$ is the ascending union of its subspaces 
$$
L_{0,n}:={\rm span} \{p^*q \; | \; p,q \in Q_n\}
$$
and each $L_{0,n}$ is a subalgebra of $L$ because $(p^*q)(x^*y)=\d_{xq}p^*y$. It is also clear that 
$L_{0,n} \subset L_{0,n+1}$ because 
$$
p^*q=\sum_{a \in s^{-1}(t(q))} p^*a^*aq
$$
(this uses the fact that $Q$ has no sinks). 

By \cite[Prop. 4.1]{XWChen2}, the linear map $kQ \to L(Q)$, $p \mapsto p$, is 
injective for any quiver $Q$. As a consequence,  there is a well-defined linear map
$$
\Phi_n:L_{0,n} \to \End_{kI} (kQ_n), \qquad \Phi_n(p^*q)(r):=rp^*q \;\;  (= \d_{rp}e_{t(p)} q),
$$
for $r\in Q_n$. Since 
$$
\Phi_n(p^*q)\Phi_n(x^*y)(r)=rx^*yp^*q = \Phi_n(x^*yp^*q)(r)
$$
$\Phi_n$ is an algebra anti-homorphism. Since $\Phi_n(p^*q)(p)=q$ and $\Phi_n(p^*q)(r)=0$ if $r \ne p$, 
$\Phi_n$ is injective. 

We will now show that $L_{0,n}$ and $ \End_{kI} (kQ_n)$
have the same dimension. This will complete the proof that $\Phi_n$ is an algebra isomorphism.
Since
\begin{equation}
\label{eq.eq.classes}
\{\hbox{non-zero } p^*q \; | \; p,q \in Q_n \}=\bigsqcup_{i \in I} \{ p^*q \; | \;p,q \in e_iQ_n\}.
\end{equation}
it follows that
$$
\dim_k L_{0,n}= \sum_{i \in I} |e_iQ_n|^2.
$$
On the other hand, 
$$
kQ_n = \bigoplus_{i \in I} e_i(kQ_n) = \bigoplus_{i \in I} ke_iQ_n
$$
and $ke_iQ_n$ is isomorphic as a left $kI$-module to a direct sum of $|e_iQ_n|$ copies of the simple $kI$-module $ke_i$. Hence
$$
\End_{kI} (kQ_n) = \bigoplus_{i \in I} \End_k( ke_iQ_n )\cong \bigoplus_{i \in I} M_{|e_iQ_n|}(k)
$$
is the direct sum of $|I|$ matrix algebras of sizes $|e_iQ_n|$, $i \in I$. This completes the proof that $L_{0,n}$ and $ \End_{kI} (kQ_n)$ have the same dimension. Hence $\Phi_n$ is an isomorphism.

Rather than counting dimensions one can give a more honest proof by observing that the elements in
$ \{ p^*q \; | \;p,q \in e_iQ_n\}$ are a set of matrix units for $\End_k( ke_iQ_n)$ with respect to the basis $e_iQ_n$. 

To complete the proof of the theorem we will show the $\Phi_n$s induce an isomorphism between the direct 
limits by showing that the diagram
$$
\UseComputerModernTips
\xymatrix{
L_{0,n} \ar[d]_{\Phi_n} \ar[rr] && L_{0,n+1}   \ar[d]^{\Phi_{n+1}}
\\
S_n \ar[rr]_{f \mapsto \id \otimes f} && S_{n+1}
}
$$
commutes. To this end, let $p^*q \in L_{0,n}$ where $p,q \in Q_n$, and let $r \in Q_n$ and $a \in Q_1$ be such that $ar \ne 0$. Thus $ar \in Q_{n+1}$ and $a \otimes r \in kQ_1 \otimes_{kI} kQ_n$. 
Going clockwise around the diagram, $\Phi_{n+1}(p^*q)(ar)=arp^*q$. Going anti-clockwise around the diagram, $(\id_{kQ_1} \otimes \Phi_{n})(p^*q)(a \otimes r)=a \otimes \Phi_n(p^*q)(r) = a \otimes rp^*q=arp^*q$. The diagram commutes.  
\end{pf}

\subsection{}
Most of the next result, but not part (6), is covered by \cite[Sect. 2.3]{AGGP}. It makes use of a construction in \cite{AGGP} that we now recall. 

Let $R$ be a ring and $\phi:(\ZZ_{\ge 1},+) \to \End_{\sRings}(R)$ a monoid homomorphism to the monoid of
ring endomorphisms of $R$.  Let $\ZZ[t_+]$ and $\ZZ[t_-]$ be the 
polynomial rings on the indeterminates $t_+$ and $t_-$. We define the ring $R[t_-,t_+;\phi]$ to be the quotient of
the free coproduct
$$
\ZZ[t_-]*_\ZZ R *_\ZZ \ZZ[t_+]
$$
modulo the ideal generated by the relations
\begin{align*}
 t_-r= & \phi(r)t_-, \quad r \in R, 
  \\
  rt_+ = & t_+\phi(r), \quad r \in R,
  \\
  t_-t_+= & \phi(1)
  \\
  t_+t_-= & 1.  
  \end{align*}

The image of the map $\phi$ in the next result is the subalgebra $t_-Lt_+=eLe$ where $e$ is the idempotent $t_-t_+$.

\begin{prop}
\label{prop.LQ+-}
We continue to assume $Q$ has neither sinks nor sources and $L$ denotes $L(Q)$. 
For each $i \in I$ pick an arrow $a_i$ ending at $i$ and define
$$
t_+:= \sum_{i \in I} a_i
\qquad \hbox{and} \qquad t_-:=t_+^*.
$$
Define a non-unital ring homomorphism $\phi:L \to L$ by $\phi(x):=t_- xt_+$. Then
\begin{enumerate}
\item{}
$t_+t_-=1$; 
  \item 
  If $n>0$, then $L_n= t_+^nL_0$ and $L_{-n}= L_0t_-^n$;
   \item 
   $L(Q)$ is generated  by $L_0$ and $t_+$ and $t_-$; 
   \item{}
   in the notation of \cite{AGGP}, $L=L_0[t_-,t_+;\phi]$;
   \item{}
   if $n$ is positive, $L_n=(L_1)^n$ and $L_{-n} =(L_{-1})^n$; 
  \item 
   $L(Q)$ is strongly graded.
\end{enumerate}
 \end{prop}
\begin{pf}
(1)
This follows from the fact that $a_ia_i^* =e_i$ and $a_ia_j^*=0$ if $i \ne j$.

(2)
Suppose $n >0$ and let $b \in L_n$ and $c \in L_{-n}$. Then $b=t_+^nt_-^nb$ and $c=ct_+^nt_-^n$, and 
$t_-^nb, ct_+^n \in L_0$. This proves (2) and (3) is an immediate consequence.

(4)
See \cite[Sect. 2 and Lem. 2.4]{AGGP}.

(5)
This is proved by induction. For example, assuming $n >0$ and starting with (2) and the induction hypothesis, $L_n=t_+^nL_0$, we have
$$
(L_1)^{n+1}=(L_1)^nL_1=t_+^nL_0L_1=t_+^nL_1 =t_+^{n+1}L_0=L_{n+1}.
$$
The proof for $L_{-n}$ is similar.

(6)
Suppose $n>0$. Then $1 =t_+^nt_-^n \in L_nL_{-n}$ so $L_nL_{-n}=L_0$.

Now 
$$
1 = \sum_{i \in I} e_i= \sum_{i \in I}  \sum_{a \in s^{-1}(i)} a^*a = \sum_{a \in Q_1} a^*a
$$
so $1 \in L_{-1}L_1$.  It now follows from (5) that $L_{-n}L_n=L_0$.
 \end{pf}

Dade's Theorem \cite[Thm. 2.8]{Dade} on strongly graded rings gives the next result.

\begin{cor}
If $Q$ has neither sinks nor sources,  there is an equivalence of categories
$$
\Gr L(Q) \equiv \Mod L_0.
$$
\end{cor}

Because $L(Q)$ is strongly graded, \cite[(2.12a)]{Dade} tells us that each $L_n$ is an invertible $L_0$-bimodule
and the multiplication in $L$ gives $L_0$-bimodule isomorphisms
$$
L_m \otimes_{L_0} L_n \stackrel{\sim}{\longrightarrow} L_{m+n}
$$
for all $m$ and $n$. In other words, if we use the multiplication in $L$ to identify $L_1^{-1}$ with $L_{-1}$ and define $L_1^{\otimes (-r)}:=(L_{-1})^{\otimes r}$ for all $r>0$, then 
$$
L= \bigoplus_{n=-\infty}^\infty L_1^{\otimes n}
$$
where the tensor product is taken over $L_0$.

 \subsection{}
 We have now completed the proof that
 $$
 \QGr(kQ) \equiv \Mod_r S(Q) \equiv \Mod_\ell L(Q)_0 \equiv \Gr_\ell L(Q)
 $$ 
 when $Q$ has no sinks or sources.
 It is possible to prove the equivalence $ \QGr(kQ) \equiv \Gr_\ell L(Q)$ directly by modifying the 
arguments in section 4 of \cite{Sm2} for the free algebra $k \langle x_0,\ldots,x_n\rangle$
so they apply to $kQ$. The required changes are minimal and straightforward so we leave the details of
the next three results to the reader.  

The next result is proved in \cite[Thm. 4.1]{AB} but the following proof is more direct. 

\begin{prop}
\label{prop.L.flat}
The ring $L$ is flat as a right $kQ$-module.
\end{prop}
\begin{pf}
Since $L$ is the ascending union of the finitely generated free right $kQ$-modules 
\begin{equation}
\label{defn.F_r}
F_n=\sum_{p \in Q_n} p^* (kQ) = \bigoplus_{p \in Q_n} p^* (kQ) 
\end{equation}
it is a flat right $kQ$-module. 
\end{pf}

A version of the following result for finitely presented not-necessarily-graded modules is 
given in \cite[Sect. 6]{AB}. 

\begin{prop}
\label{prop.L.tors}
If $M \in \Gr(kQ)$, then $L \otimes_{kQ} M=0$ if and only if $M \in \Fdim (kQ)$. 
\end{prop}
\begin{pf}
The argument in \cite[Prop. 4.3]{Sm2} works provided one replaces ``free module'' by ``projective module''.
\end{pf}

A version of the next result for finitely presented not-necessarily-graded modules is 
given in \cite[Sect. 6]{AB}.

\begin{thm}
\label{thm.Leavitt}
Let $\pi^*:\Gr(kQ) \to \QGr(kQ)$ be the quotient functor and $i^*=L \otimes_{kQ} - :\Gr(kQ) \to \Gr L$. Then 
$$
\QGr(kQ) \equiv \Gr L
$$
via a functor $\a^*:\QGr(kQ) \to \Gr L$ such that $\a^*\pi^* = i^*$.
\end{thm}
\begin{pf}
The argument in \cite[Thm. 4.4]{Sm2} works here.
\end{pf}

\subsection{}
The referee pointed out that the equivalence between $\QGr(kQ)$ and $\Gr L(Q)$ can be 
understood using ideas about perpendicular subcategories of $\Gr(kQ)$, as in \cite{GL}, and 
ideas about universal localization that are implicit in \cite{GP}.

We defined $L(Q)$ as the universal localization $\Sigma^{-1}(kQ)$. Let $\varphi:kQ \to L(Q)$ be the universal  
$\Sigma$-inverting map. As Schofield remarks, \cite[p.56]{Scho}, $\varphi$ is an epimorphism in the category of rings. It is also an epimorphism in the category of graded rings. The restriction functor $\varphi_*:\Gr(L(Q))
\to \Gr(kQ)$ therefore embeds $\Gr(L(Q))$ as a fully exact subcategory (see \cite[p. 280]{GP} for the definition) of 
$\Gr(kQ)$. 

Because $L(Q)$ is flat as a right $kQ$-module (Prop. \ref{prop.L.flat}) $\varphi:kQ \to L(Q)$ satisfies the slightly stronger property of being a homological epimorphism in the category of graded rings
(i.e., the equivalent properties of \cite[Thm. 4.4]{GL} are satisfied).


\section{Examples}
\label{sect.egs}

\subsection{} 
If  $\dim_k kQ < \infty$, then $S(Q)=0$.

\subsection{}
If  $Q$ is the cyclic quiver
$
 \UseComputerModernTips
\xymatrix{
1 \ar[r] & 2 \ar[r] & \cdots \ar[r] & n \ar@/^1.2pc/[lll]
}
$
then $S(Q) \cong k^n$.

\subsection{}
By \cite{Sm2}, if $Q$  has one vertex and $r$ arrows, then $kQ$ is the free algebra $k\langle x_1,\ldots,x_r\rangle$, the Bratteli diagram for $S(Q)$ is 
$$
\UseComputerModernTips
\xymatrix{
 1 \;\ar@<.5ex>[r]_{\cdot}   \ar@<-1ex>[r]^{\cdot}    & \; r \;  \ar@<.5ex>[r]_{\cdot}   \ar@<-1ex>[r]^{\cdot}     & \; r^2 \; \ar@<.5ex>[r]_{\cdot}   \ar@<-1ex>[r]^{\cdot}   & \; r^3 \; \ar@<.5ex>[r]_{\cdot}   \ar@<-1ex>[r]^{\cdot}     & \;  r^4 \cdots
}
$$
where there are $r$ arrows between adjacent vertices, and
$$
S(Q) \cong \liminj_n  M_r(k)^{\otimes n}.
$$

\subsection{}
Different quivers can have a common Veronese quiver. For example, 
\vskip .05in
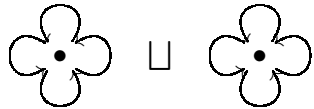
 $$
 \UseComputerModernTips
\xymatrix{
\bullet  \ar@(ul,dl)[]  \ar@(dl,dr)[]   \ar@(dr,ur)[]  \ar@(ur,ul)[]  & \bigsqcup &  \bullet \ar@(ul,dl)[]  \ar@(dl,dr)[]   \ar@(dr,ur)[]  \ar@(ur,ul)[] 
}
$$
\vskip .12in
\noindent
is the 2-Veronese quiver of both
\vskip .1in
$$
Q=\UseComputerModernTips
\xymatrix{
\bullet \ar@<-.5ex>[r]  \ar@/^1pc/[r] & \bullet \ar@<-.5ex>[l]  \ar@/^1pc/[l]
}
\qquad \hbox{and} \qquad
Q'= \quad
 \UseComputerModernTips
\xymatrix{
  \bullet  \ar@(dl,dr)[]     \ar@(ur,ul)[] 
}
\quad 
 \sqcup 
 \quad
  \UseComputerModernTips
\xymatrix{
  \bullet  \ar@(dl,dr)[]     \ar@(ur,ul)[] 
}
$$
\vskip .12in
It now follows from Theorem \ref{thm4} that $\QGr(kQ) \equiv \QGr(kQ')$ and, by the comments after 
Theorem \ref{thm4}, 
\begin{align*}
S(Q)  & = S(Q^{(2)})
\\
& = S(Q')
\\
& \cong \bigg(\liminj_n  M_2(k)^{\otimes n}\bigg) \oplus \bigg( \liminj_n  M_2(k)^{\otimes n}\bigg).
\end{align*}

\subsection{}
It is not  obvious that  $\QGr( k\langle x,y\rangle /(y^2))$ is equivalent to
$\QGr(kQ)$ for some quiver $Q$.

\begin{prop}
Let
$$
Q= \qquad \UseComputerModernTips
\xymatrix{
1    \ar@/^/[rr]  \ar@(ul,dl)[]  && 2 \ar@/^/[ll]
}
$$
The Bratteli diagram for $S(Q)$ is 
$$
\UseComputerModernTips
\xymatrix{ 
1 \ar[r] \ar[dr]  & 2 \ar[r] \ar[dr] & 3 \ar[r]  \ar[dr] &  5 \ar[r]  \ar[dr] &   8   \ar[r]  \ar[dr]&  \cdots
\\
 1 \ar[ur]  & 1  \ar[ur] & 2   \ar[ur] & 3 \ar[ur] &  5 \ar[ur] &    \cdots
}
$$
and 
$$
 \QGr(kQ) \equiv \QGr  \frac{k\langle x,y\rangle}{(y^2)}  .
$$
\end{prop}
\begin{pf}
We will use the notation $E_i$ and $E_{ij}$ that appears in the proof of Proposition \ref{prop.brat}.

The powers of the incidence matrix for $Q$ are
$$
C^n=\begin{pmatrix}
f_n      &  f_{n-1}  \\
f_{n-1}      &  f_{n-2}
\end{pmatrix}
$$
where $f_{-1}=0$, $f_0=f_1=1$, and $f_{n+1}=f_n+f_{n-1}$ for $n \ge 1$. As a $kI$-bimodule
$$
kQ_n \cong E_{11}^{f_n} \oplus    E_{12}^{f_{n-1}} \oplus  E_{21}^{f_{n-1}} \oplus  E_{22}^{f_{n-2}}
$$
and as a left $kI$-module
$$
kQ_n \cong E_1^{f_{n+1}} \oplus E_2^{f_n}.
$$
It follows that $S_n \cong M_{f_n}(k) \oplus M_{f_{n-1}}(k)$ and the Bratteli diagram is as claimed.
But this Bratteli diagram also arises in  \cite{Sm1} where it is shown that
$$
\Mod S(Q) \equiv \QGr \frac{k \langle x,y \rangle}{(y^2)}.
$$ 
It also follows from the main result in \cite{HS}, which was written after this paper, that $\QGr(k \langle x,y \rangle/(y^2))$ is equivalent to $\QGr(kQ)$.
\end{pf}

As explained in \cite{Sm1}, we can interpret $k \langle x,y \rangle/(y^2)$, and therefore $kQ$, as a non-commutative homogeneous coordinate ring of the space of Penrose tilings of the plane.  This is consistent with
Connes' view  \cite[Sect. II.3]{NCG} of the norm closure of $S(Q)$, over $\CC$, as a C$^*$-algebra coordinate ring for the space of Penrose tilings.

For each integer $r\ge 1$ there is a quiver $Q$ such that $\QGr(kQ)$ is 
equivalent $\QGr k\langle x,y \rangle/(y^{r+1})$; see section \ref{sect.multinacci} and \cite{HS}.

\subsection{}
The previous example can be generalized as follows. 

\begin{prop}
Let
$$
Q= \qquad \UseComputerModernTips
\xymatrix{
\hdot    \ar@/^/[rr]      \ar@(ul,dl)[]_m   &&  \hdot \ar@/^/[ll]
}
$$
\vskip .2in
\noindent
where there are $m$ arrows from the left-hand vertex to itself.
The Bratteli diagram for $S(Q)$ is 
$$
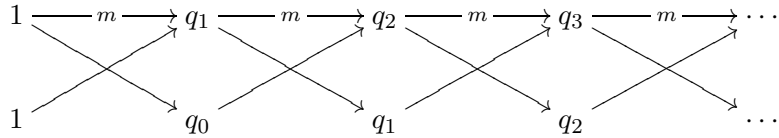
\UseComputerModernTips
\xymatrix{ 
1 \ar[rr] | {\phantom{.}m\phantom{.}} \ar[drr]  && q_1 \ar[rr] | {\phantom{.}m\phantom{.}}  \ar[drr]  && q_2 \ar[rr]
| {\phantom{.}m\phantom{.}}  
 \ar[drr]  &&   q_3 \ar[rr] | {\phantom{.}m\phantom{.}}    \ar[drr]  &&      \cdots
\\
 1 \ar[urr]  && q_0  \ar[urr] && q_1   \ar[urr] && q_2 \ar[urr] &    &    \cdots
}
$$
where there are $m$ arrows from $q_n$ to $q_{n+1}$ in the top row, and the numbers $q_n$
are given by $q_0=q_{-1}=1$ and $q_{n+1}=mq_n+q_{n-1}$ for $n \ge 0$. Furthermore, the Hilbert series of $kQ$, viewed as an element of $K_0(kI)[[t]]=(\ZZ \times \ZZ)[[t]]$, is
$$
H_{kQ}(t)=\frac{1}{1-mt+t^2}(1+t,1+(1-m)t)
$$
with the first component of $H_{kQ}(t)$ giving the multiplicity in $kQ_n$ of the simple $kI$-module that is supported at the left-most vertex. 
\end{prop}

\subsection{}
If $Q$ and $Q'$ become the same after repeatedly deleting sources and sinks, then 
$\QGr(kQ) \equiv \QGr(kQ')$ by Theorem \ref{thm.source.sink}. Therefore $S(Q)$ and $S(Q')$ are Morita equivalent, but they need not be isomorphic as the next example shows. 
The quivers $Q$ and $Q'$ are formed by adjoining a sink, respectively, a 
source, to 
\begin{equation}
\label{eg.one.loop}
Q^\circ=  \UseComputerModernTips
\xymatrix{ 
\bullet  \ar@(ur,dr)[]^x   
}
\end{equation}
By Theorem \ref{thm.source.sink}  
$$\QGr(kQ) \equiv \QGr(kQ') \equiv \QGr(kQ^\circ) \equiv \Qcoh(\Proj k[x]) \equiv \Mod k.
$$ 
In this example, $Q' = Q^{\op}$.

\begin{prop}
\label{eg.sink.source}
Let 
$$
Q=  \qquad  \UseComputerModernTips
\xymatrix{
 1  \ar[r]  \ar@(ul,dl)[]  & 2  
 }
 \qquad \hbox{and} \qquad 
 Q'= \qquad  \UseComputerModernTips
\xymatrix{
 1 \ar@(ul,dl)[]  &  \ar[l]   2.  
}
$$
 The Bratteli diagram for $S(Q)$ is 
$$
  \UseComputerModernTips
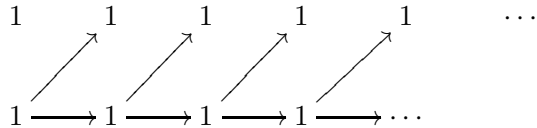
\xymatrix{
1 & 1 &1 &1 &1 &\cdots
\\
1 \ar[r] \ar[ur] & 1 \ar[r] \ar[ur] &1 \ar[r] \ar[ur] &1 \ar[r] \ar[ur] & \cdots
}
$$
and that for $S(Q')$ is
$$
  \UseComputerModernTips
\xymatrix{
2 \ar[r] & 2 \ar[r] & 2 \ar[r] & 2 \ar[r] &  \cdots
}
$$
Furthermore, $S(Q) \cong k$ and $S(Q') \cong M_2(k)$.
\end{prop}
\begin{pf}
The incidence matrices for  $Q$ and $Q'$ are 
$$
C=
\begin{pmatrix}
     1 &    0\\
      1 & 0 
\end{pmatrix} =C^n 
 \qquad \hbox{and} \qquad 
 C'=
\begin{pmatrix}
     1 &   1\\
      0 & 0 
\end{pmatrix} =(C')^n
$$
so
$$
kQ_n \cong E_{11} \oplus E_{21}
 \qquad \hbox{and} \qquad 
 kQ'_n \cong E_{11} \oplus E_{12}.
 $$
 The result follows easily from this.
 \end{pf}

Here is another way to show $\QGr(kQ) \equiv \Mod k$ for the $Q$ in Proposition \ref{eg.sink.source}. 
Write $x$ for the loop at vertex 1 and $w$ for the arrow from 1 to 2.
The path algebra is  
$$
kQ=\begin{pmatrix}
  k[x]    &  0  \\
    wk[x]  &  k
\end{pmatrix}.
$$
The two-sided ideal generated by $e_2$ is
$$
T=\begin{pmatrix}
  0   &  0  \\
    wk[x]  &  k
\end{pmatrix},
$$
which is the ideal $T$ in Lemma \ref{lem.tors.kQ},  $T$ is annihilated on the left by 
$$
\begin{pmatrix}
  k[x]   &  0  \\
    wk[x]  &  0
\end{pmatrix}
$$
so, as a left $kQ$-module, $T$ is a sum of finite dimensional modules. Therefore 
$$
\QGr(kQ)\equiv \QGr\Big( \frac{kQ}{T} \Big) \equiv  \QGr (k[x]) \equiv \Qcoh(\Proj k[x])  \equiv \Mod k.
$$

\subsection{}
Let $\sK(\cH)$ be the C$^*$-algebra of compact operators on 
an infinite dimensional separable Hilbert space  $\cH$. 
The direct limit in the category of C$^*$-algebras  of the directed 
system with Bratteli diagram (\ref{eg.K(H)+C})
is  isomorphic to $\sK(\cH) \oplus \CC \id_{\cH}$. 
The algebra $S(Q)$ in Proposition \ref{prop.K(H)+} is the 
algebraic analogue of $\sK(\cH) \oplus \CC \id_{\cH}$.

\begin{prop} 
\label{prop.K(H)+}
Let
$$
 Q= \UseComputerModernTips
\xymatrix{
{}_1 \ar@(dl,ul)[]^x & \ar[l]_w  {}_2 \ar@(dr,ur)[]_y
}
$$
The Bratteli diagram for $S(Q)$ is 
\begin{equation}
\label{eg.K(H)+C}
\UseComputerModernTips
\xymatrix{ 
1 \ar[r] & 2 \ar[r]   & 3 \ar[r] &  4 \ar[r]  &  5 \ar[r] &  \cdots
\\
1 \ar[r] \ar[ur] & 1 \ar[r] \ar[ur] & 1 \ar[r] \ar[ur] &  1 \ar[r] \ar[ur] &   1   \ar[r] \ar[ur] &  \cdots
}
\end{equation}
and 
$$
S(Q) \cong M_{\infty}(k) \oplus k.I
$$
where $M_\infty(k)$ is the algebra without unit consisting of $\NN \times \NN$ matrices with
only finitely many non-zero entries and $M_\infty(k) \oplus k.I$ is the algebra of $\NN \times \NN$ matrices that differ from a scalar multiple of the $\NN \times \NN$  identity matrix in only finitely many places.  
\end{prop}
\begin{pf}
One can see directly that 
$$ 
kQ \cong \begin{pmatrix}
 k[x]     &  k[x]\otimes w \otimes k[y]  \\
  0    &  k[y]
\end{pmatrix}.
$$
As a $kI$-bimodule, $kQ_1= kx \oplus kw \oplus ky \cong E_{11} \oplus E_{12} \oplus E_{22}$.
The $n^{\th}$ power of the incidence matrix for $Q$ is
$$
\begin{pmatrix}
 1     &    n \\
  0    &  1
\end{pmatrix}
$$
so, as a $kI$-bimodule, 
$
kQ_n \cong E_{11} \oplus E_{12}^{\oplus n} \oplus E_{22}
$
and as a left $kI$-module
$$
kQ_n \cong   E_{1}^{\oplus (n+1)} \oplus E_{2}.
$$ 
Therefore   $\End_{kI}( kQ_n)  \cong M_{n+1}(k) \oplus k$.

In order to give an explicit description of the homomorphisms 
$$\theta_n:\End_{kI}( kQ_n) \to \End_{kI}( kQ_{n+1})$$ we take
the ordered basis for $kQ_n$ consisting of the $n+2$ elements
$$
x^n, \, x^{n-1}w, \, x^{n-2}wy, \, \ldots, \, xwy^{n-2}, \, wy^{n-1}, \, y^n.
$$
The linear span of $x^n$ is a $kI$-bimodule isomorphic to $E_{11}$.
The linear span of the next $n$ elements, those with a $w$ in them, is a  $kI$-bimodule isomorphic to $E_{12}^{\oplus n}$. The linear span of $y^n$ is a $kI$-bimodule isomorphic to $E_{22}$. 
We will write an element of $f\in S_n=\End_{kI}(kQ_n)$ as 
$$
f=(A, \l) \in M_{n+1}(k) \oplus k
$$
where $A$ represents the restriction of $f$ to $E_{1}^{\oplus (n+1)}$ with respect to the ordered basis,
and $f(y^n)=\l y^n$.

The homomorphism $\theta_n:S_n \to S_{n+1}$ is defined in section \ref{sect.defn.SQ}.
In this example, 
\begin{align*}
kQ_{n+1} & =kQ_1 \otimes_{kI} kQ_n 
\\
& = \big( E_{11} \otimes_k E_1^{n+1} \big) \oplus  \big( E_{12} \otimes_k E_2  \big) \oplus  \big( E_{22} \otimes_k E_2  \big)
\\
&= x \otimes \Span\{ x^n,x^{n-1}w, \, x^{n-2}wy, \, \ldots, \, xwy^{n-2}, \, wy^{n-1} \} 
\\
& \qquad  \phantom{xxxxxxixxxx} \oplus \, (kw \otimes ky^n)   \phantom{} \oplus \, (ky \otimes ky^n) .
\end{align*}
Therefore
$$
\theta_n ( A, \l) = ( A + \l e_{n+2,n+2}, \l) = \Bigg( \begin{pmatrix} A & 0 \\ 0 & \l  \end{pmatrix}, \l \Bigg).
$$
Define $\phi_n:S_n \to M_\infty(k) \oplus kI$ by 
$$
\phi_n( A, \l): = A+\l I_{n+1}
$$
where $I_n=I-(e_{11}+ \cdots +e_{nn}) \in M_{\infty}(k) \oplus kI$.
Since $AI_{n+1}=I_{n+1}A=0$ and $I_{n+1}^2=0$, $\phi_n$ is a homomorphism of $k$-algebras sending the identity to the identity. It is straightforward to check that $\phi_{n+1}\theta_n=\phi_n$. It follows that all the
$\phi_n$s factor through a single homomorphism $\phi:S(Q) = \liminj S_n \to M_\infty(k) \oplus kI$. We 
leave the reader to check that $\phi$ is an isomorphism.
\end{pf}

\section{Relation to finite dimensional algebras with radical square zero}
\label{sect.chen}

\subsection{The work of Xiao-Wu Chen  \cite{XWChen}}
The {\sf singularity category} of a left coherent ring $R$, denoted $\sD_{\rm sg}(R)$, is the quotient of the 
derived category $\sD^b(\mod R)$ of bounded complexes of finitely presented left $R$-modules
by its full subcategory of perfect complexes. 

The following is a simplified version of the main result in \cite{XWChen}.

\begin{thm}
[X.-W. Chen]
\cite{XWChen}
Let $\L$ be a finite dimensional $k$-algebra and $J$ its Jacobson radical. Suppose $J^2=0$. 
Viewing $J$ as a left $\L$-module, define 
$$
S(\L):= \liminj \End_\L( J^{\otimes n})
$$
and the $S(\L)$-bimodule 
$$
B:=\liminj \Hom_\L(J^{\otimes n}, J^{\otimes n-1})
$$
where the maps in the directed systems are $f \mapsto \id_J \otimes f$.  Then 
\begin{itemize}
  \item 
  $B$ is an invertible $S(\L)$-bimodule with inverse $\liminj \Hom_\L(J^{\otimes n}, J^{\otimes n+1})$;
  \item 
  $S(\L)$ is a von Neumann regular ring;
  \item{}
  $J$ is a progenerator in $\sD_{\rm sg}(\L)$ with endomorphism ring $S(\L)$;
  \item 
$\Hom_{\sD_{\rm sg}(\L)}(J,-)$ is an equivalence of triangulated categories 
$$
\big(\sD_{\rm sg}(\L),[1]\big) \equiv  \big( \proj  S(\L),-\otimes_{S(\L)} B\big)
$$
where  $\proj S(\L)$ is the category of finitely generated projective right $S(\L)$-modules, and $-\otimes_{S(\L)} B$ is the translation functor on $\proj S(\L)$.
\end{itemize}
\end{thm}

If the field $k$ in Chen's theorem is algebraically closed then $\L$ is Morita equivalent to 
$kQ/kQ_{\ge 2}$ for a suitable quiver $Q$.  

\begin{thm}
\label{thm.equiv.SPS-XWC}
Let $k$ be a field, $Q$ a quiver, and $\L=kQ/kQ_{\ge 2}$. The rings $S(\L)$ and $S(Q)$ are isomorphic and 
there is an equivalence of triangulated categories 
$$
\big(\sD_{\sf sg}(\L),[1]\big) \equiv \big( \qgr(kQ), (-1) \big).
$$
\end{thm}
\begin{pf}
The Jacobson radical of $\L$ is $J=kQ_{\ge 1}/kQ_{\ge 2}$. We identify $J$ with $kQ_1$ . This identification is compatible with the $kI$-bimodule structures. 

Since $J^2=0$, 
\begin{align*}
J \otimes_{\L}   \cdots   \otimes_{\L} J = & 
J \otimes_{\L/J}  \cdots   \otimes_{\L/J} J 
\\
= & (kQ_1) \otimes_{kI}  \cdots   \otimes_{kI} (kQ_1) 
\end{align*}
so $\End_\L(J^{\otimes n}) = \End_{kI} (kQ_n)$; i.e., the individual terms in the the directed systems defining $S(\L)$ and $S(Q)$ are the same. But the maps in the directed systems are of the form $f \mapsto \id \otimes f$
in both cases so  the direct limits $S(\L)$ and $S(Q)$ are isomorphic.

Since $S(\L)$ is von Neumann regular $\proj S(\L)$ is equal to $\mod_r S(\L)$, the category of finitely presented right $S(\L)$-modules. The translation functor on $\proj S(\L)$ is $M \mapsto M \otimes_{S(\L)} B$. Hence Chen's Theorem says that
$$
\big(\sD_{\sf sg}(\L),[1]\big) \equiv  \big( \mod_r  S(Q),-\otimes_{S(Q)} B\big).
$$

We write $\Sigma$ for the translation functor $-\otimes_{S(Q)} B$. 

An auto-equivalence of an abelian category having a generator is determined its effect on the generator. The generator $\cO$ for $\qgr(kQ)$ corresponds to the generator $S(Q)$ 
under the equivalence $\Hom_{\QGr}(\cO,-):\qgr(kQ) \equiv \mod S(Q)$; since 
$\Sigma(S(Q))=B$, the auto-equivalence of $\qgr(kQ)$
that corresponds to $\Sigma$ is the unique auto-equivalence $\s$ such that 
$\Hom_{\QGr}(\cO,\s(\cO))= B$. The calculation in the next paragraph shows that $\s(\cO)=\cO(-1)$, so the auto-equivalence of $\qgr(kQ)$ that that corresponds to $\Sigma$ is $\cF \mapsto \cF(-1)$.

The equivalence $\qgr(kQ) \to \mod_r S(Q)$  sends $\cO(-1)$ to 
\begin{align*}
\Hom_{\QGr}(\cO,\cO(-1))  & = \liminj \Hom_{\Gr(kQ)}(kQ_{\ge n},kQ(-1))
\\
& = \liminj \Hom_{\Gr(kQ)}(kQ_{\ge n},kQ(-1)_{\ge n})
\\
& = \liminj \Hom_{kI}(kQ_{n},kQ(-1)_n)
\\
& = \liminj \Hom_{kI}(kQ_{n},kQ_{n-1})
\\
& = \liminj \Hom_{kI}(J^{\otimes n},J^{\otimes(n-1)})
\\
& = B.
\end{align*}
This completes the proof that $(\mod S(\L), -\otimes_{S(\L)} B) \equiv (\qgr(kQ),(-1))$.
\end{pf}

\subsection{An example}
\label{sect.multinacci}
Fix an integer $r \ge 1$,  let 
\begin{equation}
\label{Penrose.quiver}
 Q=\qquad  \UseComputerModernTips
\xymatrix{
0 \ar@(ul,dl)[] \ar[r] & 1 \ar@/^.4pc/[l] \ar[r] & 2   \ar@/^1pc/[ll]  \ar[r] & \cdots & \cdots \ar[r] & r  \ar@/^3pc/[lllll] 
}
\end{equation}
\vskip .5in
and define 
$$
\L_r := kQ/kQ_{\ge 2}.
$$

The algebra $k\langle x,y\rangle/(y^{r+1})$ in the next result is studied in \cite{Sm1}. When $r=1$ it behaves as a non-commutative homogeneous coordinate ring for the space of Penrose tilings of the plane. Thus the
equivalence of categories in the next result says that the path algebra of the quiver in (\ref{Penrose.quiver}) 
for $r=1$
is also   a  homogeneous coordinate ring for the space of Penrose tilings.

\begin{prop}
\label{prop.penrose}
The following categories are equivalent:
$$
\sD_{\sf sg}( \L_r) \equiv \qgr(kQ) \equiv \qgr  \Bigg( \frac{k\langle x,y\rangle}{(y^{r+1})} \Bigg).
$$
\end{prop}
\begin{pf}  
The incidence matrix for $Q$ is the $(r+1) \times (r+1)$ matrix
\begin{equation}
\label{M.matrix}
C= \begin{pmatrix}
  1&1&1& \cdots  & & 1   \\
    1  & 0 & 0 &   \cdots  &  & 0 \\
   0  & 1 & 0 &   \cdots    && 0 \\
   \vdots & &&&& \vdots \\
     0  & 0 & 0 &   \cdots & 1  & 0 \\
\end{pmatrix}.
\end{equation}
It follows that $kQ_1\cong E_0^{r+1} \oplus E_1 \oplus \cdots \oplus E_r$ as a left $kI$-module and, 
as a $kI$-bimodule, 
$$
kQ_1\cong E_{10} \oplus E_{21} \oplus \cdots \oplus E_{r \, r-1} \oplus \bigg( \bigoplus_{i=0}^r E_{0i} \bigg). 
$$
Thus, the dimension vector for $kQ_1$ as a left $kI$-module is $(r+1,1,\ldots,1)^\sT$ and the 
dimension vector for $kQ_n$ as a left $kI$-module is $C^{n-1}(r+1,1,\ldots,1)^\sT$.

Define $d_0=r+1$, $d_1=d_2=\cdots=d_r=1$, and $d_{n+1}=d_n+\cdots+d_{n-r}$ for $n \ge r$. 
The dimension vector for $kQ_{n+1}$ as a left $kI$-module is therefore $(d_{n+1},d_n,\ldots,d_{n-r+1})$ 
and the Bratteli diagram for $S(Q)$, written from top to bottom, is repeated copies of 
$$
  \UseComputerModernTips
\xymatrix{
S_{n-1} \ar[d]_{\theta_n} & d_{n}  \ar@{-}[dr]  \ar@{-}[d] & d_{n-1} \ar@{-}[dl] \ar@{-}[dr]  & d_{n-2}   \ar@{-}[dr]   \ar@{-}[dll]  & \cdots & d_{n-r+1} \ar@{-}[dr]  & d_{n-r}  \ar@{-}[dlllll] 
\\
S_{n} & d_{n+1} & d_{n} &  d_{n-1}   & \cdots &   & d_{n-r+1}. 
}
$$

In \cite{Sm1} it was show that 
\begin{equation}
\label{the.ring.R}
\QGr \frac{k\langle x,y\rangle}{(y^{r+1})} \equiv \Mod R_{r}
\end{equation}
where $R_{r}$ is the ultramatricial algebra associated to a Bratteli diagram that has the same underlying (unlabelled) graph as that for $S(Q)$. By Proposition \ref{prop.KG} below, $R_{r}$ and $S(Q)$ are 
Morita equivalent.
\end{pf}

 We illustrate the remark in the last paragraph of the previous proof. 
 
The Bratteli diagram for $S(\L_2)$, written from left to right,  is 
$$
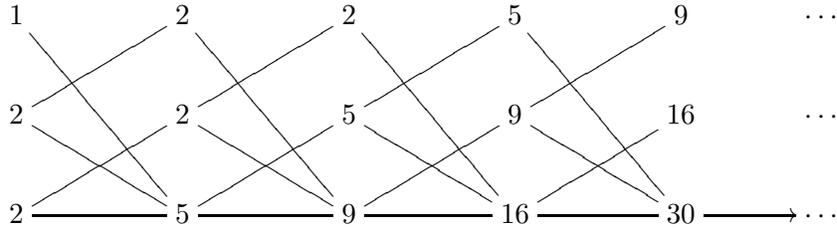
\UseComputerModernTips
\xymatrix{ 
1 \ar@{-}[ddr] & 2  \ar@{-}[ddr]  & 2 \ar@{-}[ddr]  &  5 \ar@{-}[ddr]  &  9 &  \cdots
\\
2 \ar@{-}[ur] \ar@{-}[dr] & 2 \ar@{-}[ur] \ar@{-}[dr] & 5 \ar@{-}[ur] \ar@{-}[dr]  &  9 \ar@{-}[ur] \ar@{-}[dr] & 16  &  \cdots
\\
2 \ar@{-}[r]\ar@{-}[ur] & 5  \ar@{-}[r]\ar@{-}[ur] & 9 \ar@{-}[r]\ar@{-}[ur]  &  16  \ar@{-}[r]\ar@{-}[ur] &   30   \ar[r]  &  \cdots
\\
S_1(\L_2) \ar[r] & S_2(\L_2) \ar[r] &S_3(\L_2) \ar[r] &S_4(\L_2) \ar[r] &S_5(\L_2) \ar[r] &
}
$$
The Bratteli diagram for the ring $R_2$ in (\ref{the.ring.R}), written from left to right, is 
$$
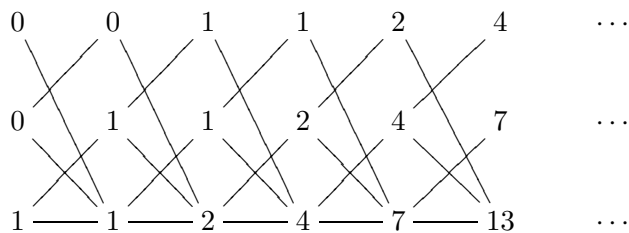
\UseComputerModernTips
\xymatrix{ 
0  \ar@{-}[ddr]  & 0  \ar@{-}[ddr] & 1   \ar@{-}[ddr] &  1 \ar@{-}[ddr]  &  2  \ar@{-}[ddr]  &  4 & \cdots
\\
0  \ar@{-}[ur] \ar@{-}[dr] & 1 \ar@{-}[ur] \ar@{-}[dr] & 1 \ar@{-}[ur] \ar@{-}[dr]  &  2 \ar@{-}[ur] \ar@{-}[dr] & 4  \ar@{-}[ur] \ar@{-}[dr] & 7 &  \cdots
\\
1 \ar@{-}[r]\ar@{-}[ur] & 1 \ar@{-}[r]\ar@{-}[ur]  & 2 \ar@{-}[r]\ar@{-}[ur] &  4  \ar@{-}[r]\ar@{-}[ur] &   7   \ar@{-}[r]\ar@{-}[ur] &   13  &  \cdots
}
$$
 
 \subsection{}
 The proof  of the next result was shown to me by
 Ken Goodearl and I thank him for allowing me to include it.
Although the result is implicit in \cite{AT}, \cite{D}, and \cite{T}, Goodearl's proof is simple and direct.

 \begin{prop}
 \label{prop.KG}
 Suppose $A$ and $B$ are ultramatricial $k$-algebras formed from Bratteli 
 diagrams on the same underlying (i.e., unlabelled) directed  graph. Then $A$ is Morita equivalent
 to $B$.
 \end{prop}
 \begin{pf}
 (Goodearl)
 The directed graph determines a directed system of free abelian groups
 whose direct limit as an ordered group is isomorphic to $K_0(A)$ and to $K_0(B)$.
Since $K_0(A)$   and $K_0(B)$ are isomorphic as ordered groups 
Elliott's results show that $A$ and $B$ are Morita equivalent. 

To see this directly, choose an ordered group isomorphism $f:K_0(A) \to K_0(B)$, 
and let $P$ be a finitely generated projective right $A$-module such that  $[P]= f^{-1}([B])$.
 Since $[B]$ is an order-unit in $K_0(B)$, $[P]$ is an order-unit in $K_0(A)$.
 This implies that $P$ is a generator in $\Mod A$. The category equivalence given by
 $ - \otimes_C P$, where $C = \End_A P$, takes the category of finitely generated
 projective right $C$-modules to the the category of finitely generated
 projective right $A$-modules with $C$ mapping to $P$. The composition
 $$
   \UseComputerModernTips
\xymatrix{
K_0(C) \ar[rr]^{-\otimes_C P} && K_0(A) \ar[rr]^f   && K_0(B)
}
$$
is an isomorphism of ordered abelian groups that sends $[C]$ to $[B]$, i.e., it is an isomorphism
of  ordered abelian groups with order unit. Elliott's   classification theorem therefore implies $C \cong B$.
But $C$ and $A$ are Morita equivalent via $P$. 
\end{pf}

\end{document}